\documentclass{amsart}
\usepackage{amssymb}
\usepackage{graphicx}
\usepackage[percent]{overpic}
\usepackage[all]{xy}
\swapnumbers
\newtheorem*{theorem*}{Theorem}
\newtheorem{theorem}{Theorem}[section]
\newtheorem{lemma}[theorem]{Lemma}
\newtheorem{proposition}[theorem]{Proposition}
\newtheorem{corollary}[theorem]{Corollary}

\theoremstyle{definition}
\newtheorem{definition}[theorem]{Definition}

\newtheorem{example}[theorem]{Example}

\theoremstyle{remark}
\newtheorem{remark}[theorem]{Remark}

\numberwithin{equation}{section}

\newcommand{\firef}[1]{Figure~{\rm\ref{#1}}}
\newcommand{\thref}[1]{Theorem~{\rm\ref{#1}}}
\newcommand{\prref}[1]{Proposition~{\rm\ref{#1}}}
\newcommand{\leref}[1]{Lemma~{\rm\ref{#1}}}
\newcommand{\coref}[1]{Corollary~{\rm\ref{#1}}}

\newcommand{\deref}[1]{Definition~{\rm\ref{#1}}}
\newcommand{\exref}[1]{Example~{\rm\ref{#1}}}

\newcommand{\reref}[1]{Remark~{\rm\ref{#1}}}

\newcommand{\seref}[1]{Section~{\rm\ref{#1}}}

\newcommand{\fig}[1]
{\raisebox{-0.5\height}%
{\includegraphics{#1}}}
\newcommand{\st}{\; | \;}                     
\newcommand{\ttt}{\otimes}                    

\newcommand{\injto}{\hookrightarrow}          
\newcommand{\isoto}{\xrightarrow{\sim}}       
\newcommand{\xxto}{\xrightarrow}              
\newcommand{\xlrto}[1]{\rlap{\mbox{$\leftarrow$}}$\xrightarrow{#1}$}

\renewcommand{\1}{\mathbf{1}}      
\renewcommand{\Vec}{\mathcal{V}ec} 
\newcommand{\GVec}{G\mathcal{V}ec} 
\newcommand{\CC}{{\mathcal C}}     
\newcommand{\CG}{{\mathcal C}/G}   
\newcommand{\Ve}{\widetilde{\mathcal V}} 
\newcommand{\V}{\mathcal V}       
\newcommand{\A}{\mathcal A}       
\newcommand{\F}{\mathcal F}       
\newcommand{\C}{\mathbb{C}}       
\newcommand{\Ctimes}{\mathbb{C}^\times}       
\newcommand{\Z}{\mathbb{Z}}       
\newcommand{\R}{\mathbb{R}}       
\newcommand{\AMod}{A{-}\text{Mod}} 
\newcommand{\VMod}{V{-}\text{Mod}} 
\newcommand{\VGMod}{V^G{-}\text{Mod}} 
\newcommand{\dd}{\mathbf{d}}       
\newcommand{\al}{\alpha}

\newcommand{\de}{\delta}

\newcommand{\la}{\lambda}

\newcommand{\ph}{\varphi}
 
\newcommand{\Rho}{R}

\newcommand{\om}{\omega}
\newcommand{\Om}{\Omega}
\newcommand{\eps}{\varepsilon}
\renewcommand{\th}{\theta}
\newcommand{\SLZ}{\mathrm{SL}_2(\Z)}

\DeclareMathOperator{\Res}{Res}
\DeclareMathOperator{\Rep}{Rep}
\DeclareMathOperator{\Ind}{Ind}
\DeclareMathOperator{\id}{id}

\DeclareMathOperator{\Hom}{Hom}

\DeclareMathOperator{\End}{End}
\DeclareMathOperator{\tr}{tr}
\DeclareMathOperator{\im}{Im}
\DeclareMathOperator{\Ker}{Ker}

\DeclareMathOperator{\Stab}{Stab}

\begin{document}

\title{On $G$--equivariant modular categories}

\author{Alexander Kirillov, Jr.}
   \address{Department of Mathematics, SUNY at Stony Brook, 
            Stony Brook, NY 11794, USA}
    \email{kirillov@math.sunysb.edu}
    \urladdr{http://www.math.sunysb.edu/\textasciitilde kirillov/}

\date{\today}

\maketitle

\section{Introduction}\label{s:intro}

This paper is was born from an attempt to provide a suitable
mathematical formalism for description of orbifold models of rational
conformal field theory. Such models arise in the study of conformal
field theories in which we have a finite group of automorphisms $G$ of
the vertex operator algebra $V$ (or, in other terminology, the chiral
algebra). In this case, we can form ``quotient'' theory which is
described by the subalgebra of invariants $V^G\subset V$. These
theories have been studied in numerous papers; see \cite{DY} for
references.

Vertex operator algebras are complicated objects, and working with
them is not easy. However, it is well known that for usual
(non-orbifold) theories, many features of the theory can be captured
by a relatively simple algebraic formalism, namely that of tensor
categories and modular functors (see, e.g.,  \cite{Tu}, \cite{BK} for
an overview).

The goal of this paper is to define a notion of a $G$-equivariant
modular tensor category that would generalize the above formalism to
the theories with a finite group of automorphisms $G$ and in
particular, give a description of the orbifold theory in terms of the
original theory. Our motivating example is the category of twisted
modules over a VOA with an action of group of automorphisms $G$.  We
plan to continue this paper with a series of papers defining
$G$-equivariant versions of modular functor, both in topological and
complex-analytic formulations.

The starting point of this paper is the definition of $G$-equivariant
fusion category due to Turaev \cite{T2} (who used the name $G$-crossed
category). This is a category with the action of a finite group $G$
and with a $G$-grading. It can be shown that for a VOA $V$ with an
action of a finite group $G$, the category of twisted $V$-modules is a
$G$-equivariant fusion category. 

The main results of this paper are as follows:

\begin{enumerate}
\item For a given $G$-equivariant fusion category $\CC$, we define the
  notion of ``orbifold'' category $\CG$, which is a fusion category
  (this definition is not new), and study basic properties of this
  category.  In the example when $\CC$ is the category of twisted
  modules over a VOA $V$, under some technical restrictions on $V$ and
  $V^G$, the orbifold category is the category of $V^G$-modules. We
  show that this construction is equivalent to the approach based on
  the notion of algebra in a category (which, in the language of
  VOA's, corresponds to considering the original VOA $V$ as a module
  over $V^G$), developed in the series of papers \cite{KO, orbi1,
    orbi2}.
\item For a $G$-equivariant fusion category $\CC$, we define the
  ``extended'' Verlinde algebra $\Ve(\CC)$ (which is no longer
  commutative) and give a simple description of the Verlinde algebra
  of the orbifold category $\V(\CG)$ in terms of $\Ve(\CC)$
  (\coref{c:verlinde7}). 
  
\item We define the $s,t$-matrices for the extended Verlinde algebra
  $\Ve(\CC)$ and show that if $s$ is non-degenerate, then these
  matrices define an action of the modular group $\SLZ$ on $\Ve(\CC)$;
  in this case, we call the category $\CC$ ``modular $G$-equivariant
  category''. Note that this definition differs from that used by
  Turaev. We show that $\CC$ is modular iff the orbifold  category
  $\CG$ is modular. 

\item   We  show that the $s$-matrix interchanges the two
  algebra structures on $\Ve$, the tensor product $\ttt$ and
  ``convolution product'' $*$. This is an analogue of the statement 
   that ``$s$-matrix diagonalizes the fusion rules'', which for
  usual fusion categories immediately gives the famous Verlinde
  formula for fusion coefficients. In the $G$-equivariant case, the
  situation is more complicated, as both $\ttt$ and $*$ are
  non-commutative, so $s$-matrix does not exactly diagonalize the
  fusion rules; however, in some special cases (e.g., when $G$ is
  commutative), one can indeed use this result to get some non-trivial
  results about the fusion coefficients of $\CC$ and $\CG$. We plan
  to pursue this in subsequent papers. 
\end{enumerate}

The paper is organized as follows. In Section 2, we define our main
object, $G$-equivariant fusion category, following \cite{T2}. In
Section 3, we explain the construction of ``orbifold quotient''
$\CG$ of a $G$-equivariant fusion category $\CC$. This orbifold
quotient is a (usual) fusion category; in the main example, when $\CC$
is the category of twisted modules over a VOA $V$, the orbifold
quotient $\CG$ is the category of modules over $V^G$, which
explains the name. In Section 4, we discuss the relation of this
construction with the approach based on algebras in category
developed in \cite{KO, orbi1, orbi2}.

Section 5 is devoted to examples
of $G$-equivariant fusion categories and corresponding fusion
categories; among other examples, it discusses category of (twisted)
$G$-graded vector spaces, whose orbifold is the category of modules
over (twisted) Drinfeld double $D^\om(G)$, and the category of twisted
modules over a VOA $V$. 

Section 6 briefly reviews presentation of morphisms in a
$G$-equivariant fusion category by ``$G$-colored'' tangles, following
\cite{T2}. 
 
In sections 7 and 8, we define extended Verlinde algebra of a $G$-equivariant
fusion category. This notion is a non-trivial generalization of the
usual Verlinde algebra; it is motivated by modular functor approach:
this generalized Verlinde algebra can be defined as a vector space
associated to torus with no punctures. Unlike the usual case, this
extended Verlinde algebra is non-commutative. 

The main results of the paper are contained in sections 9 and 10, where
we define $s$-matrix for a $G$-equivariant fusion category. Using this
$s$-matrix, we define a modular $G$-equivariant category as a
$G$-equivariant fusion category with invertible $s$. (Note that our
definition is different from Turaev's one.) We show that in any
modular $G$-equivariant category one has a natural action of the
modular group $\SLZ$ on the extended Verlinde algebra. We show that 
$\CC$ is a modular $G$-equivariant category iff the orbifold category
$\CG$ is modular; in this case, the ``untwisted sector''  $\CC_1$ is
also modular.

Throughout this paper, $G$ is a finite group.

\section{Equivariant fusion categories}\label{s:btc}
In this section, we give a definition of a $G$-equivariant tensor
category. Recall that a {\em ribbon} category is a rigid balanced
braided monoidal category $\CC$.  The words ``rigid'' and ``balanced''
mean that one has the contravariant duality functor $V\mapsto V^*$
satisfying certain properties, and that there exist functorial
isomorphisms $\de_{V}\colon V\to V^{**}$ compatible with the tensor
product; details can be found in \cite{BK}. For a braided category,
defining $\de_V$ is equivalent to defining a collection of functorial
isomorphisms $\th_V\colon V\to V$ (``twists'').

We will be mostly interested in the case when $\CC$ is a semisimple
abelian category over $\C$, with finite-dimensional spaces of
morphisms. In this case, we also assume that all functors appearing in
the definition of a ribbon category are additive and $\C$-linear on
morphisms, and that the unit object $\1$ is simple.  Such ribbon
categories are usually called {\em fusion} categories.  However, many
of our results are valid without the semisimplicity assumption.

For a fusion category $\CC$, we denote by $\V(\CC)$ the complexified
Grothendieck ring of $\CC$ (frequently also called the fusion algebra,
or the {\em Verlinde algebra}). This is a commutative algebra over
$\C$ with the basis given by classes of simple objects; it will be
discussed in detail in \seref{s:verlinde1}.

The following definition is due to Turaev \cite{T2} (who used the term
$G$-crossed category).

\begin{definition}\label{d:gfusion}
A $G$-equivariant category $\CC$ is an abelian category with the
following additional structure:

\begin{description}
\item[$G$-grading] Decomposition 
$$
\CC=\bigoplus_{g\in G} \CC_g
$$
where each $\CC_g$ is a full subcategory in $\CC$. We will call
objects $V\in \CC_g$ ``$g$-twisted''. In particular, objects $V\in
\CC_1$ will be called ``neutral''; in physical literature, the
subcategory $\CC_1$ is usually called the ``untwisted sector''

\item[Action of $G$] For each $g\in G$, we are given a  functor
  $\Rho_g\colon \CC\to \CC$ and functorial isomorphisms
  $\al_{gh}\colon \Rho_g\circ \Rho_h\isoto \Rho_{gh}$ such that
  $\Rho_1=\id$, $\Rho_g \CC_h\subset \CC_{ghg^{-1}}$, and $\al_{g_1g_2,
    g_3}\circ \al_{g_1,g_2} =\al_{g_1, g_2g_3}\circ \al_{g_2, g_3} $
  (both sides are functorial isomorphisms
  $\Rho_{g_1}\Rho_{g_2}\Rho_{g_3}\isoto \Rho_{g_1g_2g_3}$).
\end{description} 
Following \cite{T2}, we will also frequently use notation ${}^gV$ for
$\Rho_g(V)$.

A $G$-equivariant fusion category is a semisimple $G$-equivariant
abelian category which in addition has the following structures:
\begin{itemize}
\item Structure of a rigid monoidal category such that\\
  $\1$ is a simple object\\
   $\Rho_g$ is a tensor functor\\
   for $X\in \CC_g, Y\in \CC_h, X\otimes Y\in \CC_{gh}$

\item Functorial isomorphisms $\de_{V}\colon V\to V^{**}$, satisfying
  the same compatibility conditions as in the absence of $G$ (see
  \cite{BK}) and additional condition
  $\Rho_g(\de_V)=\de_{\Rho_g(V)}$. 

\item A collection of functorial isomorphisms $R_{V,W}\colon V\ttt
  W\to {}^gW\ttt V$ for every $V\in \CC_g, W\in \CC_h$, satisfying
  an analog of the pentagon axiom (see \cite[Section 2.2]{T2}). 
\end{itemize}
\end{definition}
The definition immediately implies that $\1\in \CC_1$ and that for
$V\in \CC_g$, $V^*\in \CC_{g^{-1}}$. Also, since in a rigid category
the unit object and dual is unique up to unique isomorphism, we have
canonical identification 
\begin{equation}\label{e:chig}
\begin{aligned}
&{}^g\1=\1,\\ 
&({}^gV)^*={}^g(V^*)
\end{aligned}
\end{equation}
\begin{remark}
  From now one, we will refer to associativity, unit, and $\de$
  morphism, as well as isomorphisms \eqref{e:chig}, and their
  compositions, as ``canonical'' isomorphisms, and we will omit them
  in the formulas, writing, e.g., $V_1\ttt V_2\ttt V_3$ rather than
  $(V_1\ttt V_2)\ttt V_3$. Thus, all identities between morphisms
  written below only make sense after insertion of appropriate
  ``canonical'' morphisms. Pedantic reader may complete all
  computations by inserting appropriate canonical morphisms. See
  \cite[Section 1.1]{BK} for discussion of this. 
\end{remark}

Note that any abelian category $\CC$ admits a trivial $G$-grading,
with $\CC_1=\CC, \CC_g=0$ for $g\ne 1$. Thus, as a special case this
definition includes fusion categories with action of $G$. 

As usual, existence of morphism $\de_V\colon V\to V^{**}$ is
equivalent to a system of twists $\th_V$.

\begin{lemma}\label{l:twist}
Let $\CC$ be a $G$--equivariant fusion category. Then one can define a
collection of functorial morphisms $\th_V\colon V\to {}^gV, V\in
\CC_g$, satisfying the following conditions:
\begin{enumerate}
\item $\th_\1=\id$
\item $\th_{U\ttt V}=(\th\ttt \th)R_{{}^gV, U}R_{U,V} $
\item $\th_{V^*}=\Rho_{g^{-1}}(\th_V^*)$
\item $\th_{{}^hV}=\Rho_h(\th_V)$ 
\end{enumerate}

Conversely, $\de_V$ can be recovered from $\th$, $R$, and monoidal
structure. 
\end{lemma}
The proof is completely parallel to the one in $G=\{1\}$ case (see,
e.g., \cite[Section 2.2]{BK}).

{}From now on, $\CC$ will denote a $G$--equivariant fusion category. 

The following lemma is an immediate consequence of the definition.
\begin{lemma}\label{l:fusion2}
Let $\CC$ be a $G$-equivariant fusion category. Then the fusion
algebra $\V(\CC)$ has a natural structure of a module over the Drinfeld
double $D(G)$. 
\end{lemma}

A number of examples of $G$-equivariant fusion categories is given in
\seref{s:examples}. The most important of them is the example of
twisted modules over a vertex operator algebra. 

Note that a $G$--equivariant fusion category is {\em not} a fusion
category: the braiding $R$ does not satisfy the usual axioms of the
commutativity isomorphism. However, for any $G$--equivariant fusion
category $\CC$, there are two related fusion categories. The first one
is the ``untwisted sector''  $\CC_1\subset \CC$  (in
terminology of \cite{T2}, {\em neutral category}): it is easy to see
from the definition that the equivariant fusion structure on $\CC$,
when restricted to $\CC_1$, defines a (usual) fusion category
structure on $\CC_1$.

The second way to construct a a fusion category from a
$G$--equivariant fusion category is by ``orbifolding'', or taking, in
appropriate sense, the quotient by the action of $G$. This
construction is studied in detail in the next section.

\section{Orbifold category}\label{s:orbicat}
Let $\CC$ be a $G$-equivariant fusion category. In this section we
define the notion of ``orbifold'' category. This construction had
appeared implicitly in \cite{DY} and was explicitly defined in
\cite{orbi2}, to which the reader is referred for proofs (where it is
denoted by $\CC^G$); in the case of ``free'' action of $G$, it is also
described in \cite{CM}.

\begin{definition}\label{d:orbicategory} 
Let $\CC$ be a $G$-equivariant fusion category. Then the
``orbifold fusion category'' $\CG$ is a category with the following
objects and morphisms: 

\begin{itemize}
\item Objects: pairs $(X, \{\ph\}_g)$, where $X\in \CC$ and $\ph_g$ is
  a collection of $\CC$-morphisms $\ph_g \colon {}^g X\simeq X$ such
  that 
\begin{equation}\label{e:ph}
\begin{aligned}
&\ph_1=\id\\
&\ph_g \Rho_g(\ph_h)=\ph_{gh}
\end{aligned}
\end{equation}

\item Morphisms: $(X,\{\ph_g\})\to (Y, \{\psi_h\})$ are $\CC$-morphisms
  $f\colon X\to Y$ such that $\psi_g\circ R_g(f)=f \circ\ph_g$. 
\end{itemize}
\end{definition}

\begin{remark}\label{r:ph}
Condition \eqref{e:ph} ensures that there is a canonical way of
identifying all twists ${}^gX$. Namely, if we denote, for $g,h\in G$,
\begin{equation}\label{e:ph2}
\ph_{g,h}=\ph_g^{-1}\ph_h\colon {}^hX\isoto {}^gX,
\end{equation}
then \eqref{e:ph} implies 
\begin{equation}\label{e:ph3}
\begin{aligned}
&\ph_{g,h}\ph_{h,f}=\ph_{g,f}\\
&\Rho_g\ph_{a,b}=\ph_{ga,gb}
\end{aligned}
\end{equation}
Conversely, if $\ph_{g,h}\colon {}^hX\isoto {}^gX$ is a system of
isomorphisms satisfying \eqref{e:ph3}, then 
$\ph_g=\ph_{1,g}$ satisfies \eqref{e:ph}. 
\end{remark}

\begin{lemma}
The orbifold category $\CG$ is an abelian category.
\end{lemma}
The proof of this lemma is straightforward (see \cite{orbi2} for details).

It is immediate from the definition that $\CG$ has a natural
structure of a module category over the category $\Rep G$ of
finite-dimensional complex $G$-modules: for a representation
$\rho\colon G\to \End V$, we define $\rho\ttt (X,\{\ph\})=(V\ttt X,
\{\rho_g\ttt\ph_g\})$ (see \cite{O} for overview of the notion of
module category).

The  notation $\CG$  is justified by the following lemma:

\begin{lemma}
  If $\CC$ is semisimple, and $G$ acts freely on the set $I(\CC)$ of
  isomorphism classes of simple objects, then $\CG$ is also
  semisimple, with the set of simple objects $I(\CG)=I(\CC)/G$, and
  fusion algebra $\V(\CG)=\V(\CC)^G$. 
\end{lemma}

However, in general it is not true that $\V(\CG)=\V(\CC)^G$;
instead, $\V(\CG)$ can be described in terms of {\em extended}
Verlinde algebra of $\CC$ (see \coref{c:verlinde7}). 

The following theorem, proved in \cite{DY} in the context of twisted
modules over VOAs and in \cite{orbi2} in the language of
$G$-equivariant categories, gives a full description of $\CG$ as an
abelian category.

For $i\in I(\CC)$, let $G_i=\Stab(i)=\{g\in G\st
{}^gV_i\simeq V_i\}$. Note, however, that a priori we do not have a
canonical isomorphism between ${}^gV_i$ and $V_i$. We can choose such an
isomorphism $\la_{g}\colon {}^gV_i\simeq V_i $ arbitrarily; then we
have $\la_g \Rho_g(\la_h)=\al_{gh}\la_{gh}$ for some $\al\colon
G_i^2\to \Ctimes$. It is easy to see that $\al$ is a two-cocycle, and
that the class $[\al]\in H^2(G_i, \Ctimes)$ does not depend on the
choice  of $\la_g$. Thus, action of $G$ by automorphisms of $\CC$ not
only gives an action of $G$ on the set $I(\CC)$, but also, for every
$i\in I(\CC)$, a cohomology class $[\al_i]\in H^2(G_i,\Ctimes)$. Such a
cohomology class defines a central extension 
$$
1\to \Ctimes\to \widehat G^\al_i\to G_i\to 1
$$
Define the twisted group algebra $\C^{\al_i}[G_i]$ by 
\begin{equation}\label{e:Galpha}
\C^{\al_i}[G_i]=\C[\widehat G^\al_i]/([c]-c[1])
\end{equation} 
where, for $c\in \Ctimes$, we denote by $[c]$ the class of the
corresponding element in the central extension $\Ctimes\subset
\widehat G^\al_i$. Then one easily sees that choosing a lifting
$G_i\injto \widehat G^\al_i$ gives a basis in $\C^\al[G]$ consisting
of of classes $[g], g\in G$ with multiplication law
$[g][h]=\al_{gh}[gh]$.
Thus, representations of $\C^{\al_i}[G_i]$ (or, equivalently, the category of
representations of the central extension $\widehat G^\al_i$ such that
an element $c\in \Ctimes\subset \widehat G^\al_i$ acts by
multiplication by $c$) are exactly the projective representations of
$G_i$ with the cocycle $\al$.

\begin{theorem}[\cite{DY, orbi2}]\label{t:abelian}
As an abelian category, $\CG$ is equivalent to 
$$
\bigoplus_{i\in I/G} \Rep \C^{\al_i^{-1}}[G_i]
$$
where $i$ runs over the set of representatives of $G$-orbits in $I$. 
\end{theorem}
Proof and details can be found in \cite[Theorem 3.5]{DY},
\cite[Theorem 3.5]{orbi2}. 

\begin{corollary}\label{c:semisimple}
$\CG$ is a semisimple abelian category. 
\end{corollary}
Indeed, it is known (see, e.g., \cite{Kar}) that $\C^{\al_i}[G_i]$ is
a semisimple associative algebra, and thus the category of
representations is semisimple.

\begin{corollary}\label{c:verlinde}
As a vector space, the Verlinde algebra $\V=\V(\CG)$ is given by  
$$
\V=\bigoplus_{i\in I/G}Z(\C^{\al_i^{-1}}[G_i])
$$
where 
$$
Z(\C^{\al_i^{-1}}[G_i])=(\C^{\al_i^{-1}}[G_i])^{G_i}=\{f\in
\C^{\al_i^{-1}}[G_i] \st
[g]f[g]^{-1}=f\ \forall g\in G_i\}
$$
is the center of $\C^{\al_i^{-1}}[G_i]$
\end{corollary}
\begin{proof}
Since $\C^{\al_i^{-1}}[G_i]$ is semisimple,
$Z(\C^{\al_i^{-1}}[G_i])=\bigoplus_{\la} \End (V_\la)$, where $V_\la$
are irreducible $\C^{\al_i^{-1}}[G_i]$-modules. Thus,
$\id_{V_\la}\mapsto [V_\la]$ is an isomorphism between
$Z(\C^{\al_i^{-1}}[G_i])$ and the Verlinde algebra of
$\Rep \C^{\al_i^{-1}}[G_i]$. 
\end{proof}
\begin{remark}
In fact, the isomorphism of \coref{c:verlinde} is actually an
isomorphism of associative algebras if $\V$ is considered as an
algebra with respect to the convolution product, which will be
discussed in detail in \seref{s:verlinde1}.
\end{remark}

It turns out that the structure of $G$-equivariant fusion category on
$\CC$ gives rise to structure of fusion category on $\CG$.

\begin{theorem}\label{t:orbicategory2}
Let $\CC$ be a $G$-equivariant fusion  category. Then the orbifold
category $\CG$ is itself a fusion category, with the following
tensor product, unit object, and duality:  
\begin{align*}
&(X,\{\ph_g\})\ttt (Y, \{\psi_g\})=(X\ttt Y, \{\ph_g\ttt\psi_g\})\\
&\1_{\CG}=(\1, \{\id\})\\
&(X,\{\ph_g\})^*=(X^*,\{(\ph^*_g)^{-1}\})
\end{align*}
Equivalently, in terms of morphisms $\ph_{g,h}$ defined in
\reref{r:ph}, the corresponding isomorphisms for  $X^*$ are defined by
$\ph_{h,g}^*\colon {}^hX^*\isoto {}^gX^*$.  

The associativity, unit, and balancing isomorphisms $\de$ are
inherited from $\CC$ and the braiding is defined by
$$
X\ttt Y\xxto{R}{}^g Y\ttt X\xxto{\psi_g\ttt 1}Y\ttt X.
$$

In this category, the universal twist $\th\colon X\to X$ is defined as
follows: if we write $X=\bigoplus_h X_h$, $X_h\in \CC_h$, then $\th$ is 
the direct sum of the following composition:   
$$
X_h\xxto{\th}{}^hX_h\xxto{\ph_h}X_h. 
$$

\end{theorem}
The proof is straightforward: the associativity, unit, and rigidity
automorphisms are inherited from $\CC$; verification of all identities
is left to the reader.

One has some natural functors relating categories $\CC$ and $\CG$: 
\begin{equation}\label{e:functors}
\begin{aligned}[t]
\Ind\colon \CC&\to \CG\\
           V&\mapsto (X,\{\ph_g\}\\
           X&=\bigoplus_h{}^h V\\
           \ph_g&\colon \bigoplus_h{}^{gh}V\to \bigoplus_h{}^{h}V
\end{aligned}\qquad 
\begin{aligned}[t]
         \Res\colon \CG&\to \CC\\
           (X,\{\ph_g\})&\mapsto X
\end{aligned}
\end{equation}
where $\ph_g\colon \bigoplus_h{}^{gh}V\to \bigoplus_h{}^{h}V$ is the
permutation of summands.

\begin{example}\label{x:ind1}
Let $V=\1\in \CC$; then $\Ind V=\bigoplus_g {}^g\1$. Using canonical
isomorphism \eqref{e:chig} to identify ${}^g\1= \1$, we can identify  
\begin{equation}\label{e:ind1}
\Ind V=\bigoplus_g {}^g\1=\F(G)\ttt \1,
\end{equation}
where $\F(G)$ is the vector space of functions on $G$. One easily
sees that under this identification, $\ph_g$ becomes the left
regular action of $G$ on $\F(G)$.
\end{example}

 The
functor $\Ind$ should be thought of as an analog of ``induction''
functor in usual representation theory; the functor $\Res$, as the
``restriction'' functor (see \exref{x:NtimesG} below).

The following theorem is an analogue of \cite[Theorem 1.6]{KO}. 

\begin{theorem}\label{t:functors}         
\begin{enumerate}
\item The functors $\Ind$ ans $\Res$ are adjoint to each other: for
  any $V\in \CC, X\in \CG$. 
\begin{equation}
\Hom_\CC(\Res X, V)=\Hom_{\CG}(X, \Ind V)
\end{equation}
\item The functor $\Res$ is a tensor functor.
\item $\Res$ is compatible with duality and balancing: $(\Res
  V)^*=\Res (V^*)$, $\Res(\de)=\de$. In particular, $\dim_{\CC}\Res X=
  \dim_{\CG} X$ 
\end{enumerate}
\end{theorem}
The proof is straightforward and left to the reader.

\section{Orbifolds and algebras in a category}\label{s:algebras}

The theory of orbifold categories is closely related to the approach
based on algebras in category, discussed in \cite{orbi1, orbi2,
  fuchs}. For reader's convenience, we briefly review here some of the
definitions and results of \cite{KO, orbi1, orbi2}. Recall that if
$\A$ is a fusion category, a commutative algebra in $\A$ is an object
$A\in\A$ with multiplication morphism $\mu\colon A\ttt A\to A$
satisfying obvious axioms. We will always assume that $A$ satisfies
two additional properties:
\begin{equation}\label{e:algebra1}
\begin{aligned}
&\th_A=\id\\
&A\text{  is rigid }
\end{aligned}
\end{equation}
The last condition means that $\1$ has multiplicity 1 in $A$ and the
composition\\
 $A\ttt A\xxto{\mu} A\to \1$ is a non-degenerate pairing
(see \cite{KO} for details).

For each such algebra $A$, we can define the category of (left)
$A$-modules, with objects being objects of $\A$ with a morphism
$\mu_V\colon A\ttt V\to V$. We denote the category of $A$-modules by
$\AMod$. It can be shown that $\AMod$ is a semisimple abelian category
(see \cite[Theorem 3.3]{KO}).  This category has a structure of a
monoidal category which in general is not braided. One also has
natural functors $F\colon \A\to \AMod, G\colon \AMod\to \A$ (see
\cite[Theorem 1.6]{KO}) which are adjoint to each other; the functor
$F$ is a tensor functor.

Let us consider a special case, when we are given an action of a
finite group $G$ by automorphisms $\pi_g$ of $A$. Assume additionally
that this action satisfies the following conditions:
\begin{equation}\label{e:algebra2}
\begin{aligned}
&\text{The action is faithful: if }g\ne 1,\text{ then } \pi_g\ne \id \\
&A^G=\1.
\end{aligned}
\end{equation}

\begin{example}\label{x:F(G)}
  Let $\A=\Rep G$ be the category of $G$-modules, and $A=\F(G)$ the
  algebra of functions on $G$, with pointwise multiplication and
  structure of $G$ module given by left regular action of $G$ on
  $\F(G)$. Let $\pi_g$ be the {\em right} regular action of $G$ on
  $\F(G)$; it commutes with the left regular action and thus defines
  an algebra automorphism. One easily sees that this action satisfies
  properties  (\ref{e:algebra1}, \ref{e:algebra2}). 
\end{example}

For an algebra with action of $G$ by automorphisms, we can define the
notion of $g$-twisted module: an $A$-module is called $g$-twisted if
$\mu_V\circ R^2=\mu_V\circ (\pi_{g^{-1}}\ttt \id))$.

 The following
theorem summarizes many of the results on \cite{orbi1, orbi2}.

\begin{theorem}\label{t:algebra1}
Let $\A$ be a fusion category and $A$---a commutative algebra in
 $\A$ satisfying conditions {\rm(\ref{e:algebra1})}, with an action of a
 finite group $G$ satisfying conditions  {\rm(\ref{e:algebra2})}. Then:
\begin{enumerate}
\item \cite[Theorems 2.11, 2.15]{orbi1} $\A$ contains as a full
  subcategory the category $\Rep G$. The algebra $A$ lies in this
  subcategory and can be identified with the algebra $\F(G)$ of
  \exref{x:F(G)}.  
  
\item \cite[Section 5]{orbi2} The category $\CC=\AMod$
  has a natural structure of a $G$-equivariant fusion category, with
  $\CC_g$ being the subcategory of $g$-twisted $A$-modules.
  
\item \cite[Theorems 4.1, 4.4]{orbi2} The category $\A$ is naturally
  equivalent to the orbifold $(\AMod)/G$. Under this equivalence, the
  functors $F\colon \A\to \AMod$, $G\colon \AMod\to \A$ are identified
  with the functors $\Res, \Ind$ defined by \eqref{e:functors}. 
\end{enumerate}
\end{theorem}

It turns out that this result can be reversed.
\begin{theorem}\label{t:algebra2}
  Let $\CC$ be a $G$-equivariant fusion category,  $\A=\CG$; by
  \thref{t:orbicategory2}, $\A$ is a fusion category. Let 
  $A=\Ind(\1_\CC)\in \A$. Then $A$ has a natural structure of a
  commutative algebra in $\A$ with an action of $G$, satisfying
  conditions  {\rm(\ref{e:algebra1}), (\ref{e:algebra2})}, and the
  category of $A$-modules is naturally equivalent to the category
  $\CC$.
\end{theorem}
\begin{proof}
  By \exref{x:ind1}, as an object of $\CG$, $A=\F(G)\otimes
  \1_{\CG}$, where $\F(G)$ is the algebra of functions on $G$.
  Define multiplication on $A$ and action of $G$ by automorphisms on
  $A$ as in \exref{x:F(G)}. Then it is immediate that $A$ is a
  commutative algebra in $\A$, satisfying conditions
  (\ref{e:algebra1}), (\ref{e:algebra2}).
  
  One can also explicitly describe the category of $A$-modules.
  Namely, writing $A=\bigoplus {}^g\1$, we see that an $A$-module in
  the category $\A=\CG$ is an object $(X,\{\ph_g\})\in \CG$
  along with a collection of $\CC$-morphisms $\mu_g\colon X\to X$ such
  that $\mu_{g_1}\mu_{g_2}=\de_{g_1,g_2}\mu_{g_1}, \sum \mu_g=\id$ and 
$$
\ph_g\Rho_g\mu_h =\mu_{gh}\ph_g.
$$

To show that the category of $A$-modules is equivalent to $\CC$,
construct the functors $\CC\to \AMod, \AMod\to \CC$ as follows:

$\CC\to\AMod$: for $V\in \CC$, consider $\Ind V=(\bigoplus {}^g
V,\{\ph_g\})\in \CG$ and define on it the structure of $A$-module
by defining $\mu_g=\id_{{}^gV}$. It is easy to check that it satisfies
all the required properties. 

$\AMod\to \CC$: let $(X, \{\ph_g\})\in \CG$, with structure of
$A$-module given by $\mu_g\colon X\to X$. Define an object $V\in \CC$
by $V=\im \mu_1$. 

It is trivial to check that these two functors are inverse to each
other, and thus establish an equivalence of categories $\CC\simeq
\AMod$.

\end{proof} 

Thus, starting from a fusion category $\A$ with a commutative algebra
$A$ satisfying (\ref{e:algebra1}), (\ref{e:algebra2}), one can
construct a $G$--equivariant fusion category $\CC=\AMod$.  The
original category $\A$ can be recovered as $\A=\CG$.  Conversely,
starting from a $G$-equivariant fusion category $\CC$, one can define
the fusion category $\A=\CG$ and an algebra $A\in \A$ with an
action of $G$; then $\CC$ can be recovered as the category of
$A$-modules.  This allows us to use some results about the category of
$A$-modules to answer questions about relation between $\CC$ and
$\CG$.

\section{Examples of equivariant fusion categories}\label{s:examples}

\begin{example}\label{x:vec}
  Let $\CC=\Vec$ be the category of vector spaces, with trivial
  grading and trivial action of $G$: $\Rho_g=\id$ for all $g$. Then
  the corresponding orbifold category $\Vec/G$ is the category $\Rep
  G$ of finite-dimensional $G$-modules, with the usual tensor product.
  The functors $\Ind$ and $\Res$ are given by $\Ind(V)=V\ttt \F(G)$,
  where $\F(G)$ is the space of functions on $G$ with left regular
  action of $G$, and $\Res(X)=X$ is the forgetful functor.
\end{example}

\begin{example}\label{x:vecG}
  Let $\CC=\GVec$ be the category of $G$-graded vector spaces. It can
  be explicitly described as the category with simple objects
  $X_g,g\in G,$ and tensor product, duality given by $X_g\ttt
  X_h=X_{gh}, X_g^*=X_{g^{-1}}$. Define the action of $G$ by $\Rho_g
    X_h=X_{ghg^{-1}}$. Then orbifold category $\GVec/G$ is the
    category of finite-dimensional modules over the Drinfeld double
    $D(G)$ with the usual tensor product. In this case, $\Res\colon
  \Rep D(G)\to \GVec$ is the
  usual forgetful functor, and $\Ind\colon \GVec\to \Rep D(G)$ is the
  usual induction functor $\Vec\to \Rep G$ with suitably defined
  $G$-grading. This example is discussed in slightly different
  language in \cite[Sections 5,6]{orbi1}.  
\end{example}

\begin{example}\label{x:twistedvecG}
  Let $\CC$ be a twisted category of $G$-graded vector spaces, i.e. a
  rigid monoidal category which coincides with $\GVec$ as an abelian
  category, and has tensor product, duality defined so that $X_g\ttt
  X_h\simeq X_{gh}$, $ X_g^*\simeq X_{g^{-1}}$ (non-canonically). Each
  such category defines a 3-cocycle $\om\in C^3(G, \Ctimes)$: if we
  choose a system of isomorphisms $\al_{gh}\colon X_g\ttt X_h\isoto
  X_{gh}$, then $\al_{g_1,g_2g_3}\al_{g_2,g_3}=\om(g_1, g_2,
  g_3)\al_{g_1g_2, g_3}\al_{g_1,g_2}$. It is easy to show that two
  such categories are equivalent as monoidal categories iff
  $[\om]=[\om']$, and that this defines a bijection between
  equivalence classes of twisted categories of $G$-graded vector
  spaces and $H^3(G, \Ctimes)$. We will denote by $\GVec^\om$ the
  twisted category of $G$-graded vector spaces with cocycle $\om$.  (A
  slightly different but equivalent description of these categories is
  given in \cite[Section 1.3, Section 2.6]{T2}).
    
  Define the action of $G$ by $\Rho_g X= X_g\ttt X \ttt X_g^*$, and
  the braiding isomorphism as composition
$$
X_g\ttt X_h\isoto X_g\ttt X_h \ttt X_g^*\ttt X_g ={}^gX_h\ttt X_g.
$$
One easily sees that this defines on $\GVec^\om$ a structure of
$G$-equivariant fusion category. In this case, the corresponding
orbifold category $\GVec^\om/G$ can be shown to coincide with the category
of modules over twisted Drinfeld double $D^\om (G)$ as defined in
\cite{DPR}, \cite{DN}. 

  Note that unlike the construction in \cite[Section 2.6]{T2}, our
  definition does not require that the cocycle $\om$ be
  $G$-invariant. On the other hand, our construction is not the most
  general: e.g., there are different ways to define braiding in
  $\GVec^\om$; see \cite{T2} for details. 
\end{example}

\begin{example}\label{x:NtimesG}
  Let $N$ be a group on which $G$ acts by automorphisms. Let $\CC=\Rep
  N$ be the category of finite-dimensional $N$-modules. Define on this
  category an action of $G$ as follows: for a module $M$, we let
  $R_g(M)$ be the same vector space as $M$ but with the action of $M$
  defined by $\rho_{R_g(M)}(n)=\rho_{M}(g^{-1}(n))$. Then $\CC$
  becomes an equivariant fusion category (with trivial $G$-grading),
  and the corresponding orbifold category is $\Rep N/G=(G\ltimes
  N)-Mod$ (see \cite[Theorem 2.1]{orbi2}). In this case, the functors
  $\Res\colon \Rep (G\ltimes N) \to \Rep N$ and $\Ind\colon \Rep N\to
  \Rep (G\ltimes N)$ are the usual induction and restriction functors.
\end{example}

\begin{example}\label{x:twistedVOA}
  Let $V$ be a rational vertex operator algebra such that the category
  of $V$-modules is a fusion category, and let $G$ be a finite group
  subgroup of automorphisms of $V$. Then we can define the category
  $\CC=\VMod_{tw}$ of twisted $V$-modules as in \cite{DVVV}, or, in
  more detail, in \cite{DLM1}. This category by definition has a
  $G$-grading: $\VMod_{tw}=\bigoplus_g \VMod_g$, with $\VMod_1=\VMod$.
  For a twisted $V$-module $M$, let $R_g(M)$ coincide with $M$ as a
  vector space but define the action of $V$ by
  $Y_{R_g(M)}(v,z)=Y_M(g^{-1}(v),z)$.  Then the category $\VMod_{tw}$
  has a natural structure of a $G$-equivariant category.  Under
  suitable assumptions on $V$, this category has a natural structure
  of $G$-equivariant fusion category. Indeed, as shown in \cite{KO,
    orbi1}, the category of twisted modules over $V$ is equivalent to
  the category of $A$-modules, where $A$ is $V$ considered as an
  associative commutative algebra in the category $\A=\VGMod$ of
  moduleks over $V^G$. Thus, applying the results of \seref{s:algebras}, we
  see that $\VMod_{tw}$ is a $G$-equivariant fusion category. As a
  corollary of \thref{t:algebra1}, we see that the orbifold category
  in this case is $\A=\CG$, i.e., $\VGMod=\VMod/G$.
  
  This example is the main motivation for the study of
  $G$-equivariant fusion categories.  
\end{example} 
\begin{remark}
  It would be interesting to give a direct description of the monoidal
  structure on the category of twisted modules over $V$, i.e. a
  description which does not use restriction of the modules to $V^G$.
  The only result in this direction we were able to find is the paper
  \cite{fusion}, so many details are still missing. 
\end{remark}

\section{Graphical description of morphisms in $G$-equivariant
  categories} \label{s:tangles}

In this section, we briefly review the graphical technique for
representing morphisms in a $G$-equivariant fusion category,
generalizing well-known graphical technique for representing morphisms
in a braided tensor category by tangles. The results of this section
are due to Turaev \cite{T2}. For now, we present a very simplified
description; a more detailed description, using the language of links
and tangles in a 3-manifold with a principal $G$-bundle, can be found
in \cite{T2}. 

Recall that a {\em tangle diagram} is a collection of oriented arcs
and circles in $\R\times[0,1]$, where the arcs ends must lie on the
lines $\R\times\{0\}, \R\times\{1\}$. The only intersections allowed
are transversal double intersections, and for each such intersection,
one of the strands is specified as ``top'', and the other as
``bottom''. Such tangle diagrams naturally arise as projections of
tangles in $\R^3$, and it is well-known that two diagrams correspond
to isotopic tangles iff they can be obtained one from another by a
sequence of Reidemeister moves. 

This can be generalized to $G$-equivariant situation. For a tangle
diagram $T$, a {\em segment} of $T$ is a part of an arc or a circle
between two undercrossings (i.e., points where a given arc goes under
one of the other arcs).

\begin{figure}[tb]
\begin{overpic}
{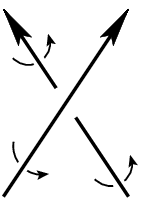}
\put(25,12){$g$}
\put(70,20){$h$}
\put(-30,53){$ghg^{-1}$}
\put(70,-10){$W$}
\put(5,-10){$V$}
\put(-15,100){${}^gW$}
\end{overpic}
\hspace{2cm}
\begin{overpic}
{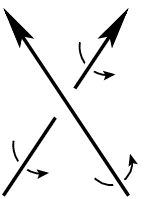}
\put(25,12){$g$}
\put(70,20){$h$}
\put(60,60){$h^{-1}gh$}
\put(70,-10){$W$}
\put(5,-10){$V$}
\put(65,100){${}^{h^{-1}}V$}
\end{overpic}
\caption{Crossing}\label{f:crossing}
\end{figure}

\begin{definition}\label{d:color}
Let $\CC$ be a $G$-equivariant category, and $T$ --- a tangle
diagram. A $\CC$ coloring of $T$ is an assignment to every segment of 
an arc or circle a pair $(g, V), g\in G, V\in \CC_g$ (a {\em color}) 
satisfying the two conditions below. In the figures, we will show a 
color by writing the object $V$ next to the segment, and writing $g$ 
next to an arrow going under the corresponding strand, as in 
\firef{f:crossing}. Since $V$ determines $g$, we will frequently write 
just $V$ and omit notation of $g$.  

\begin{enumerate}

\item For every circle, the ordered product $\prod_i g_i=1$, where the
  product is over all undercrossings of this circle, and $g_i$ are the
  colors of crossing strand.

\item For two segments separated by an undercrossing, the colors are
related as shown in \firef{f:crossing}. 

\end{enumerate}
\end{definition}

\begin{remark}
  Condition (1) ensures compatibility of condition (2). Indeed,
  condition (2) implies that with each undercrossing, $V$ is replaced
  by ${}^{g_i}V$. Thus, if we have a circle, condition (2) implies
  that  ${}^gV=V$, where $g=\prod_i g_i$.
  However, in a tensor category, condition ${}^gV=V$ is almost never
  satisfied; at best, one could expect that ${}^gV$ is isomorphic to
  $V$, but then we would need to specify explicitly the choice of
  isomorphism, since there is no canonical isomorphism. To avoid these
  problems, we require that $\prod g_i=1$. More general situation,
  with $g\ne 1$ and a choice of isomorphism ${}^gV\isoto V$, is best
  dealt with by allowing circles with coupons described below.
\end{remark}

It can be shown (see \cite[Lemma 3.2.1]{T2}) that a coloring is
uniquely determined by specifying the color of just one segment on
every arc and circle of $T$ (these colors can not be chosen
arbitrarily: condition (2) imposes restrictions on them). Therefore,
in many figures  we will specify the color at just one point on each
circle or arc. 

We have an obvious action of $G$ on the set of coloring of a given
tangle diagram $T$: for a color $(V,h)$, we define action of $g$ by
$R_g(V,h)=({}^gV, ghg^{-1})$. One easily checks that this preserves
conditions (1), (2).   

As in the non-equivariant case, we can assign morphisms to tangle
diagrams as follows. Let $T$ be a tangle diagram. Then bottom of $T$
defines a sequence of triples $(g,V,\eps)$, where $(g,V)$ is the color
of the corresponding segment of $T$, and $\eps=\pm$ is defined by the
direction of the corresponding segment: $\eps=+$ if is is directed up,
$\eps=-$ if it is directed down. Define 
\begin{equation}
X_{in}(T)=\bigotimes V_i^{\eps_i},
\end{equation}
where the tensor product is over all ends of the arcs at the bottom of
$T$, in the natural order (left to right), and $V^\eps=V$ for
$\eps=+$, $V^\eps=V^*$ for $\eps=-$. If the bottom of $T$ is empty,
we let $X_{in}(T)=\1$. In a similar way, we define
$X_{out}(T)$ by taking the product over the top of the diagram $T$.

\begin{figure}[ht]
\begin{overpic}
{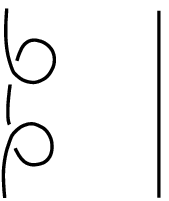}
\put(30,45){\xlrto{\Om_1}}
\end{overpic}
\hspace{1cm}
\begin{overpic}
{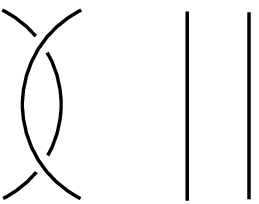}
\put(40,35){\xlrto{\Om_2}}
\end{overpic}
\hspace{1cm}
\begin{overpic}
{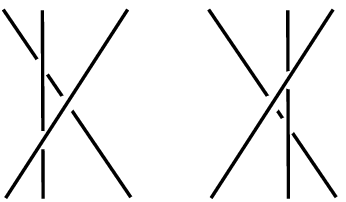}
\put(40,25){\xlrto{\Om_3}}
\end{overpic}
\caption{Reidemeister moves}\label{f:redem}
\end{figure}

\begin{theorem}\label{t:tangle}
  Let $\CC$ be a $G$-equivariant fusion category, Then there is a
  unique way to assign to every colored tangle diagram $T$ a
  $\CC$-morphism $F(T)\colon X_{in}(T)\to X_{out}(T)$ so that the
  following properties are satisfied:
\begin{enumerate}
\item For elementary crossing, ``cap'' and ``cup'' diagrams, $F(T)$ is
  the commutativity morphism $R$, rigidity morphism $V\otimes V^*\to \1$,
  and $\1\to V\ttt V^*$ respectively.

\item $F(T_1\ttt T_2)=F(T_1)\ttt F(T_2)$, and $F(T_1\circ
  T_2)=F(T_1)F(T_2)$, where $T_1\ttt T_2$ is the diagram obtained by
  placing $T_2$ to the right of $T_1$, and $T_1\circ T_2$ is the
  diagram obtained by placing $T_1$ on top of $T_2$ {\rm(}it is  defined
  only if $X_{out}(T_2)=X_{in}(T_1)$\rm{)}. 
\end{enumerate}

So defined assignment $F$ satisfies the following properties:
\begin{description}
\item[$G$-equivariance] For any $g\in G$,
$$
F(\Rho_g T)=\Rho_g(F(T))\colon {}^gX_{in}(T)\to {}^g X_{out}(T)
$$
\item[Independence of the choice of projection]
 $F(T)$ is invariant under the Reidemeister moves for framed tangles
 shown in \firef{f:redem}.
 
\item[Orientation reversal] $F(T)$ is invariant under simultaneously
  reversing direction of a component of $T$ and replacing its color
  $(V,g)$ by $(V^*, g^{-1})$.

\item[Tensor product]
Replacing in $T$ a component with color $(V_1\ttt V_2, g_1g_2)$ by
two components, obtained by doubling the original component, and  with
colors $(V_1, g_1)$ and $(V_2, g_2)$, does not change $F(T)$.  

\end{description}

\end{theorem}

The condition of invariance under Reidemeister moves means that $F(T)$
depends only on the isotopy class of the tangle and not on the choice
of diagram representing this tangle. However, accurate explanation of
this requires that we introduce the notion of $\CC$-colored tangle
which takes some time; see \cite{T2}.

\begin{example}
For tangle $T$ shown in \firef{f:twist} with $V\in \CC_g$, $F(T)=\th_V\colon
V\to {}^gV$ (which explains why $\th$ is called ``twist''). For
technical reasons, however, we will usually just draw box labeled by
$\th$ instead of drawing the twist. 

\begin{figure}[ht]
\fig{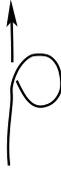}
\caption{Graphical representation of the universal twist}\label{f:twist}
\end{figure}
\end{example}

\begin{figure}[ht]
\raisebox{-0.5\height}{\begin{overpic}
{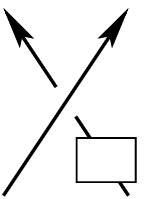}
\put(70,-10){$W_1$}
\put(43,35){$W_2$}
\put(49,13){$\Phi$}
\put(5,-10){$V$}
\put(-15,100){${}^gW_2$}
\end{overpic}}\quad
=
\quad
\raisebox{-0.5\height}{\begin{overpic}
{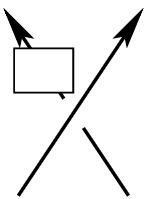}
\put(70,-10){$W_1$}
\put(13,62){$ {}^g\Phi$}
\put(11,-10){$V$}
\put(-15,100){${}^gW_2$}
\end{overpic}}
\qquad 
\raisebox{-0.5\height}{\begin{overpic}
{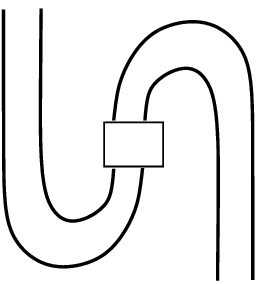}
\put(41,45){$\Phi$}
\end{overpic}}
\quad=
\quad 
\raisebox{-0.5\height}{\begin{overpic}
{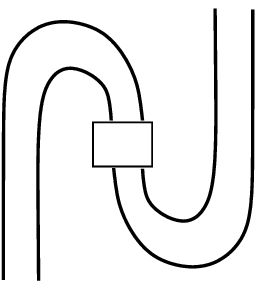}
\put(38,45){$\Phi$}
\end{overpic}}
\caption{}\label{f:coupon}
\end{figure}

These results can naturally be extended to graphs with coupons, i.e.
rectangular boxes; the color of such a coupon should be a morphism
$\Phi\colon X_{in}\to X_{out}$, where $X_{in}, X_{out}$ are tensor
products of colors at the bottom (respectively, top) of the coupon.
So defined invariants for graphs with coupons satisfy, in addition to
the Reidemeister moves above, the moves shown in \firef{f:coupon}. See
\cite{T2} for proofs and details.

\section{Verlinde algebra 1}\label{s:verlinde1} 
In this section we recall known facts about the Verlinde algebra of a
fusion category $\CC$. All of them are well-known; however, they are
formulated here in somewhat unusual form for the convenience of later
generalization to the $G$-equivariant case, which will be done in the
next section.

Throughout this section,  $\CC$ is  a fusion category. For simplicity,
we assume that the set $I$ of isomorphism classes of simple objects in
$\CC$ is finite (all the results are actually valid without this
assumption if we allow some of the objects we consider to be infinite
sums of simple objects in $\CC$). 

\subsection*{Definition of $\V$}
We define complex vector space
$\V=\V(\CC)$ by
\begin{equation}
\V=\bigoplus_{i\in I} \Hom(V_i,V_i),
\end{equation}
where the sum is over the set $I=I(\CC)$ of equivalence classes of simple
objects in $\CC$. Then $\V$ has a natural basis $\chi_i=[\id_{V_i}]$.

This vector space also has a more invariant definition.

\begin{theorem}\label{t:verlinde0}
  $\V$ is isomorphic to the vector space spanned by classes $[\ph]$,
  where $\ph\colon V\to V, V\in \CC$, with the following relations
\begin{enumerate}
\item For any $\la\in \C,
\ph,\psi\in \End V$, one has 
\begin{equation}\label{e:verrel1}
\la[\ph]=[\la\ph],\quad [\ph+\psi]=[\ph]+[\psi].
\end{equation}
\item For any $\ph\in\End V$ and isomorphism $f\colon V\isoto V'$, one
  has 

\begin{equation}\label{e:verrel2}
[f\ph f^{-1}]=[\ph]
\end{equation}
\item If $W=\bigoplus W_i$, for some $W_i\in \CC$, and $\ph\colon W\to
  W$, $\ph=\sum \ph_{ij}, \ph_{ij}\colon W_i\to W_j$, then 
\begin{equation}\label{e:verrel3}
[\ph]=\sum [\ph_{ii}]
\end{equation}

\end{enumerate}
\end{theorem}
\begin{proof}
  Denote temporarily by $\V'$ the vector space generated by $[\ph]$
  with relations \eqref{e:verrel1}--\eqref{e:verrel3}. Define a map
  $\V'\to \V$ as follows. Let $\ph\colon V\to V$; write $V=\bigoplus
  H_i\ttt V_i$, where $H_i=\Hom(V_i, V)$ are multiplicity spaces. Then
  $\ph=\bigoplus \ph_i\ttt\id_{V_i}$, where $\ph_i\colon H_i\to H_i$
  is a linear map. We define the map $\V'\to \V$ by 
\begin{equation}
[\ph]\mapsto \sum (\tr\ph_i) \chi_i. 
\end{equation}
It is easy to see that relations \eqref{e:verrel1}--\eqref{e:verrel3}
are satisfied, so this gives a well-defined map $\V'\to \V$. This map
is clearly surjective. It is also surjective: if $\tr \ph_i=0$ for all
$i$, then it follows from \eqref{e:verrel3} that $[\ph]=0$ in $\V'$. 
Thus, this map gives an isomorphism $\V'\isoto \V$. 
\end{proof}

{}From now on, we will frequently use this theorem, writing various
operations in $\V$ in terms of classes $[\ph]$. Of course, whenever we
define something in terms of classes $[\ph]$, we need to verify that
relations \eqref{e:verrel1}--\eqref{e:verrel3}. This is usually
trivial and therefore we will not write it explicitly. 

In particular, for an object $V$ we define $[V]=[\id_V]=\sum (\dim
H_i)\chi_i$. This brings us back to the standard definition of $\V$ as
the complexified Grothendieck ring of $\CC$.

\subsection*{Tensor product in $\V$}
The vector space $\V$ has a structure of associative commutative
algebra, which we will denote by $\ttt$. It is defined by
\begin{equation}
[\ph]\otimes [\psi]=[\ph\otimes \psi]
\end{equation}

In particular, 
\begin{equation*}
\chi_i\ttt \chi_j=\sum N_{ij}^k\chi_k
\end{equation*}
where $N_{ij}^k$ are fusion coefficients: $V_i\ttt V_j\simeq \sum
N_{ij}^kV_k$. 

The unit with respect to $\ttt$ is $\chi_0=[\1]$.

\subsection*{Convolution product in $\V$} 
There is also another associative, commutative product in $\V$, which
we will denote by $*$ (it is sometimes called ``the convolution
product''). It is defined as follows: for $\ph\colon V_i\to V_i,
\psi\colon V_j\to V_j$, we define
\begin{equation}
[\ph]*[\psi]=\de_{ij}d_i^{-1}[\ph \psi]
\end{equation}
where 
\begin{equation}\label{d:di}
d_i=\dim V_i.
\end{equation}
Recall that the numbers $d_i$ are always non-zero (see,
e.g. \cite[Section 2.4]{BK}).

In particular, 
\begin{equation}
\chi_i*\chi_j=d_i^{-1}\de_{ij}\chi_i.
\end{equation}

The unit with respect to $*$ is
\begin{equation}\label{e:d}
\dd=\sum_{i\in I} d_i\chi_i.
\end{equation}

\subsection*{Bilinear form}
The Verlinde algebra $\V$ also has a natural bilinear form. Namely,
for $\ph\colon V\to V$, let $\ph^*\colon V^*\to V^*$ be the adjoint
morphism. This defines an algebra automorphism $V\to V$.  Also, define
the ``constant term'' of an element $x\in \V$ by $[\chi_i]_0=0, i\ne
0$ and $[\chi_0]_0=1$. Then we define the bilinear form on $\V$ by
\begin{equation}\label{e:bilinform}
(\ph,\psi)=[\ph\ttt\psi^*]_0. 
\end{equation}

\begin{lemma}\label{l:bilinform}
The bilinear form \eqref{e:bilinform} has the following properties:
\begin{enumerate}
\item It is symmetric and non-degenerate
\item $(\chi_i, \chi_j)=\de_{ij}$
\item  $(x\ttt y,z)=(x,z\ttt y^*).$
\end{enumerate}
\end{lemma}
The proof of this lemma is straightforward and left to the reader. 

\subsection*{Dimension homomorphism}
We define the dimension homomorphism $d\colon \V\to \C$ by
letting, for $\ph\colon V\to V$,
\begin{equation}\label{e:dim1}
d([\ph])=tr_V(\ph).  
\end{equation}
Thus, $d[V]=\dim V$, so $d(\chi_i)=\dim(V_i)=d_i$.  One easily sees
that the dimension map has the following properties:
\begin{equation}\label{e:dim2}
\begin{aligned}
&d(x\ttt y)=d(x)d(y),\\
&d(x^*)=d(x),\\
&d(\1)=1.
\end{aligned}
\end{equation}

\begin{lemma}\label{l:dim}
  Let $\dd \in \V$ be defined by \eqref{e:d}. Then 

\begin{enumerate}
\item $(\dd,x)=d(x)$
\item $\dd^*=\dd$.
\item For any $x\in \V$, $\dd \ttt x=d(x)\dd $. 
\end{enumerate}
\end{lemma}
\begin{proof}
(1) is immediate from \leref{l:bilinform} and definition of $\dd$; (2)
    follows from (1) and \eqref{e:dim2}. To prove (3), note that it
    suffices that both sides have the same inner product with any
    $y\in \V$. Using results of \leref{l:bilinform} and \eqref{e:dim2},
    we can write 
$$
(\dd\ttt x,y)=(\dd, y\ttt x^*)=d(y\ttt x^*)=d(y)d(x)=(d(x)\dd, y).
$$
\end{proof}

\begin{lemma}
For $x,y\in \V$, one has 
$$
(x,y)=d(x*y). 
$$
\end{lemma}
\begin{proof}
Suffices to prove this for $x=\chi_i, y=\chi_j$, in which case
$d(\chi_i*\chi_j)=\de_{ij}d_i^{-1}d(\chi_i)=\de_{ij}$. 
\end{proof}

\begin{corollary}
  If $V$ is a simple module, and $\ph, \psi\in
\End(V)$, then 
\begin{equation}
([\ph],[\psi])=\frac{1}{\dim V}\tr (\ph\psi).
\end{equation} 
\end{corollary}
Note that this may be false if $V$ is not simple: for example, if
$\CC$ is the category of vector spaces, then $([\ph],[\psi])=(\tr
\ph)(\tr \psi)$, which in general is not equal to $\frac{1}{\dim V}\tr
(\ph\psi)$  

Notice that there is a nice symmetry between $\ttt$ and $*$:
\begin{alignat*}{3}
&\1\ttt x=x  \qquad  & & \dd\ttt x=d(x)\dd  \qquad  & & (x,y)=[x\ttt y^*]_0\\
&\1*x=[x]_0\1\qquad  & & \dd *x=x \qquad  & & (x,y)=d(x*y)
\end{alignat*}
This symmetry interchanges $\1$ with $\dd$ and $[\ ]_0$ with $d$.

\section{Verlinde algebra 2}\label{s:verlinde2} 
In this section we give a definition of the extended Verlinde algebra
for $G$-equivariant fusion categories, following the same steps as in
\seref{s:verlinde1} with suitable changes. The results of this section
are new.

\subsection*{Definition of $\Ve$}
\begin{definition}\label{d:verlinde}
  Let $\CC$ be a $G$-equivariant fusion category, and $V_i$
  ---representatives of isomorphism classes of simple objects in
  $\CC$. Then the extended Verlinde algebra of $\CC$ is defined by 
$$
\Ve(\CC)=\bigoplus_{i\in I, g\in G}\Hom_{\CC}(V_i, {}^g V_i).
$$
\end{definition}

This definition is motivated by modular functor point of view:
$\Ve(\CC)$ is the vector space assigned to a torus with no boundary
components (see \cite[Section 8.6]{T2}).

{}From now on, we assume that the set $I$ of isomorphism classes of
simple objects in $\CC$ is finite, so that $\Ve$ is finite-dimensional.

Note that if $V_i\in \CC_h$, then $\Hom(V_i, {}^gV_i)=0$ unless
$ghg^{-1}=h$. Thus, the Verlinde algebra can be written as follows:
\begin{equation}
\begin{aligned}
&\Ve(\CC)=\bigoplus_{g,h:gh=hg} \Ve_{g,h}(\CC)\\
&\Ve_{g,h}(\CC)=\bigoplus_{i\in I_h}\Hom(V_i, {}^g V_i)
\end{aligned}
\end{equation}
where $I_h$ is the set of isomorphism classes of simple objects in
$\CC_h$. In particular, 
$$
\Ve_{1,*}=\bigoplus_{i\in I} \Hom(V_i, V_i)=\V(\CC)
$$
is the usual Verlinde algebra, i.e. the complexified Grothendieck
ring of the category $\CC$.
\begin{remark}
If the action of the group $G$ on the set of isomorphism classes of
simple objects is free, then it immediately follows form the
definition that $\Ve_{g,*}=0$ for $g\ne 1$, and thus $\Ve=\V$,
i.e. the extended Verlinde algebra coincides with the usual Verlinde
algebra. 
\end{remark}

As before, $\Ve$ has a more invariant definition. 
\begin{theorem}\label{t:verlinde20}
  Let $g,h\in G, gh=hg $. Then $\Ve_{g,h}$ is isomorphic to the vector space
  spanned by classes $[\ph]$, where $\ph\colon V\to {}^gV, V\in \CC_h$, with
  the following relations
\begin{enumerate}
\item For any $\la\in \C,
\ph,\psi\colon V\to {}^gV$, one has 
\begin{equation}\label{e:verrel21}
\la[\ph]=[\la\ph],\quad [\ph+\psi]=[\ph]+[\psi]`.
\end{equation}
\item For any $\ph\colon V\to {}^gV$ and isomorphism $f\colon V\isoto V'$, one
  has 

\begin{equation}\label{e:verrel22}
[R_g(f)\,\ph f^{-1}]=[\ph].
\end{equation}
\item If $W=\bigoplus W_i$, for some $W_i\in \CC_h$, and $\ph\colon W\to
  W$, $\ph=\sum \ph_{ij}, \ph_{ij}\colon W_i\to {}^gW_j$, then 
\begin{equation}\label{e:verrel23}
[\ph]=\sum [\ph_{ii}].
\end{equation}

\end{enumerate}
\end{theorem}
\begin{proof}
The proof is parallel to the proof of \thref{t:verlinde0}, with the
following modification: the map $\Ve'\to \Ve$ is defined by 
$$
[\ph]\mapsto \sum (\tr \ph_{ii})[\la_{ii}]
$$
where $\ph\colon V\to {}^gV$, $V=\bigoplus H_i\ttt V_i$,
$H_i=\Hom(V_i, V)$, and we write $\ph=\bigoplus \ph_{ij}\ttt
\la_{ij}$, $\ph_{ij}\colon H_i\to H_j, \la_{ij}\colon V_i\to {}^g
V_j$.
\end{proof}
{}From now on, we will frequently use this theorem, writing various
operations in $\Ve$ in terms of classes $[\ph]$. Of course, whenever
we define something in terms of classes $[\ph]$, we need to verify
that relations \eqref{e:verrel21}--\eqref{e:verrel23} are satisfied.
This is usually trivial and therefore we will not write it explicitly.

In particular, for any $V\in \CC$ we define $[V]=[\id_V]\in \Ve$ and
denote $\chi_i=[V_i]\in \Ve_{1,*}$. The elements $\chi_i$ form a basis
in $\Ve_{1,*}$.

We have an obvious action of $G$ on $\Ve$, given by 
\begin{equation}\label{e:Rx}
\begin{aligned}
\Rho_x\colon \Ve_{g,h}&\to \Ve_{xgx^{-1}, xhx^{-1}}\\
[\ph]&\mapsto [\Rho_x\ph].
\end{aligned} 
\end{equation}

\subsection*{Tensor product in $\Ve$}

The vector space $\Ve$ has a structure of associative algebra, which
we will denote by $\ttt$. It is defined as follows: if $[\ph]\in
\Ve_{g_1,h_1}, [\psi]\in\Ve_{g_2,h_2}$, then
\begin{equation}\label{e:product1}
[\ph]\ttt [\psi]=\begin{cases}
[\ph\ttt\psi]\in \Ve_{g_1, h_1h_2}& g_1=g_2\\
0& g_1\ne g_2
\end{cases}
\end{equation}

\begin{lemma}
\begin{enumerate} 
\item The product $\ttt$ defined by \eqref{e:product1} defines on
  $\Ve$ a structure of associative algebra with unit 
\begin{equation}\label{e:unit2}
\tilde\1 =\sum\chi^g_0,
\end{equation}
where $\chi_0^g\colon \1\to {}^g\1$ is the canonical isomorphism
\eqref{e:chig}. 

\item For $\ph\in \Ve_{g_1,h_1},
\psi\in\Ve_{g_2,h_2}$, one has 
\begin{equation}\label{e:commutativity}
[\ph]\ttt[\psi]=[\Rho_{h_1}\psi]\ttt[\ph]. 
\end{equation}
In particular, $[\ph]\ttt [\psi]=[\psi]\ttt[\ph]$ if one of $[\ph],
[\psi]$ is in $\Ve_{*,1}$.

\item For each $g\in G$, $R_g$ is an algebra automorphism  with
  respect to $\ttt$:\\  $R_g(x\ttt y)=R_g(x)\ttt
  R_g(y)$, $R_g(\tilde \1)=\tilde \1$. 
\end{enumerate}
\end{lemma}
\begin{proof}
The proof is straightforward and is left to the reader.
\end{proof}

\subsection*{Convolution product}
Vector space $\Ve$ also has another associative product, which we will
denote by $*$ (the convolution product). It is defined as follows:

\begin{definition}\label{d:product2}
Let $\ph\colon V_i\to {}^{g_1} V_i$, $\psi\colon V_j\to
{}^{g_2}V_j$. Then $[\ph]*[\psi]\in \Ve$ is defined by 
\begin{itemize}
\item 
 If ${}^{g_2}V_j, V_i$  are not isomorphic, then  $[\ph]*[\psi]=0$
\item  If there exists an isomorphism $\la\colon {}^{g_2}V_j\isoto
 V_i$, then 
$$
[\ph]*[\psi]=d_i^{-1}
[V_j\xxto{\psi}{}^{g_2}V_j\xxto{\la}
  V_i\xxto{\ph} {}^{g_1}V_i\xxto{\Rho_{g_1}(\la^{-1})} {}^{g_1g_2}V_j ]
$$
where, as before, $d_i=\dim V_i=\dim {}^{g_2}V_j=d_j$. 
\end{itemize}
\end{definition}

It is obvious that $[\ph]*[\psi]$ is independent of the choice of
isomorphism $\la\colon {}^{g_2}V_j\isoto V_i$ and that, for $[\ph] \in
\Ve_{g_1, h_1}, [\psi]\in \Ve_{g_2,h_2}$, we have
$$
\begin{alignedat}{2}
&[\ph]*[\psi]\in \Ve_{g_1g_2,h_1}\qquad&&h_1=h_2\\
&[\ph]*[\psi]=0&&h_1\ne h_2.
\end{alignedat}
$$

This definition is chosen so that if $V$ is a simple object,
$\psi\colon V\to {}^gV$, $\ph\colon {}^gV\to {}^{hg}V$, then 
$$
[\ph]*[\psi]=\frac{1}{\dim V} [\ph\psi].
$$

\begin{lemma}
\begin{enumerate} 
\item The product $*$ defined by \deref{d:product2} defines on
  $\Ve$ a structure of an associative algebra with unit 
\begin{equation}\label{e:d2}
\dd=\sum_{i\in I}d_i\chi_i.
\end{equation}

\item  $[\ph]*[\psi]=[\psi]*[\ph]$ if one of $[\ph],
[\psi]$ is in $\Ve_{1,*}$.
\item For each $g\in G$, $R_g$ is an algebra automorphism  with
  respect to $*$:  $R_g(x*y)=R_g(x)*
  R_g(y)$, $R_g(\dd)=\dd$.
\end{enumerate}
\end{lemma}
The proof of this lemma is trivial and left to the reader. 

It is possible to give an explicit description of the structure of
$\Ve$ as an algebra with respect to $*$. Recall that for $i\in I$, we
denoted $G_i=\Stab(i)$ and that action of $G$ on $\CC$ defines for
every $i$ a cohomology class $[\al]\in H^2(G_i, \Ctimes)$ (see
discussion preceding \thref{t:abelian}).
\begin{theorem}\label{t:algebra3}
One has a canonical isomorphism of associative algebras 
$$
\Ve=\bigoplus_{i\in I} \C^{\al^{-1}}[G_i]
$$
where $\Ve$ is considered with respect to $*$ product, and
$\C^{\al^{-1}}[G_i]$ is the twisted group algebra \eqref{e:Galpha}.
\end{theorem}
\begin{proof}
  For every $i\in I, g\in G_i$, choose an isomorphism $\la_{i,g}\colon
  {}^gV_i\isoto {V_i}$. Then $\la_g R_g(\la_h)=\al_{gh}\la_{gh}$,
  where $\al\in C^2(G_i, \Ctimes)$ is the two-cocycle from
  \thref{t:abelian}. 

  Rewrite this in the following form:
$$
(R_g
(\la_h^{-1})\la_g^{-1}\la_h)(\la_h^{-1})=\al_{gh}^{-1}\la_{gh}^{-1}. 
$$
Denoting $x_g=d_i[\la_g^{-1}]\in \Ve_{g,*}$ and using
\eqref{e:verrel22}, we get $x_g*x_h=\al_{gh}^{-1}x_{gh}$.
\end{proof}

\begin{example}
Let $\CC=\GVec$ be the category of $G$-graded vector spaces as in
\exref{x:vecG}. Then, for any pair $g,h\in G$ such that $gh=hg$, we
have a canonical isomorphism $\chi_{g,h}\colon X_h\to
{}^gX_h$. Elements $\chi_{g,h}, gh=hg$ form a basis in  $\Ve$. In
this basis, the multiplication is given by
\begin{align*}
\chi_{g_1,h}*\chi_{g_2,h}&=\chi_{g_1g_2,h}\\
\chi_{g,h_1}\ttt \chi_{g,h_2}&=\chi_{g,h_1h_2}
\end{align*}
(all other products are zero). 
\end{example}

\subsection*{Bilinear form}
We  define the bilinear form on $\Ve$ as follows. First, recall that
for $\ph\colon V\to {}^gV$, we have an adjoint morphism $\ph^*\colon
{}^gV^*\to V^*$. This defines on $\Ve$ a linear map $*\colon
\Ve_{g,h}\to \Ve_{g^{-1}, h^{-1}}$ such that $(\ph\ttt
  \psi)^*=\psi^*\ttt\ph^*$ and $(\ph*\psi)^*=\psi^**\ph^*$.

Next, define the constant term map $[\ ]_0\colon \Ve\to\C$ as follows
\begin{align*}
&[\ph]_0=0, \quad \ph\colon V_i\to {}^gV_i, \ i\ne 0\\
&[\chi^g_0]_0=1,
\end{align*}
where, as before, $\chi_0^g\colon \1\to {}^g\1$ is the canonical
isomorphism. One easily sees that it completely determines $[\ ]_0$
and that $[x]_0=0$ if $x\in \Ve_{g,h}, h\ne 1$.

Now define the bilinear form 
\begin{equation}\label{e:bilinform2}
(\ph,\psi)=[\ph\ttt \psi^*]_0.
\end{equation}
\begin{lemma}\label{l:bilinform3}
The bilinear form \eqref{e:bilinform2} has the following properties: 
\begin{enumerate}
\item 
For $\ph\in \Ve_{g_1,h_1}, \psi\in \Ve_{g_2,h_2}$, we have 
$(\ph,\psi)=0$ unless $g_1=g_2^{-1}, h_1=h_2$.
\item The form is symmetric: $(\ph,\psi)=(\psi,\ph)$, non-degenerate,
  and $G$-invariant. 
\item $(\chi_i,\chi_j)=\de_{ij}$
\item $(x\ttt y,z)=(x,z\ttt y^*)$. 
\end{enumerate}
\end{lemma}
The proof of this lemma is straightforward and left to the
reader. Note, however, that $(\chi_i,\chi_j)=\de_{ij}$ is not
sufficient to determine $(\ , )$. 

\begin{example}
  Consider the subalgebra in $\Ve$ generated by classes $\chi_0^g$.
  Then this subalgebra, considered with $*$ product is isomorphic to
  the group algebra $\C[G]$, and considered with the $\ttt$ product,
  it is isomorphic to the algebra $\F(G)$ of functions on $G$. The
  bilinear form $(\ , \ )$ restricted to this subalgebra coincides
  with the standard bilinear form on $\C[G]$:
$$
(\chi_0^g, \chi_0^h)=\de_{g,h}.
$$
\end{example}

\subsection*{Dimension homomorphism}
We define the dimension homomorphism $d\colon
\Ve\to \C$ as follows: for $\ph\in \Ve_{g,h}$, we let 
\begin{equation}\label{e:dim3}
d(\ph)=\begin{cases}
\tr \ph, & g=1\\
0 &g\ne 1
\end{cases}
\end{equation}
As before, we have $d_i=d(\chi_i)=\dim V_i, i\in I(\CC)$. It is
immediate from the definition that so defined dimension homomorphism
satisfies properties similar to those in $G=\{1\}$ case with
additional property of being $G$-invariant:

\begin{equation}\label{e:dim4}
\begin{aligned}
&d(x\ttt y)=d(x)d(y),\\
&d(x^*)=d(x),\\
&d(\1)=1,\\
&d(\Rho_g\ph)=d(\ph).
\end{aligned}
\end{equation}

\begin{lemma}\label{l:dim6}
Let $\dd\in \Ve_{1,*}$ be defined by \eqref{e:d2}. Then 
\begin{enumerate}
\item $(\dd, \ph)=d(\ph)$. 
\item $\dd^*=\dd$
\item For any $x\in \Ve_{1,*}$, $\dd\ttt x=x\ttt \dd=d(x)\dd$.
\item $(\dd_h)^*=\dd_{h^{-1}}$. 
\item Let $\dd_h=\sum_{i\in I_h}d_i\chi_i \in \Ve_{1,h}$, so that
  $\dd=\sum_h \dd_h$. Then for any $x\in \Ve_{1,h}$, 
$$
x\ttt \dd_{h^{-1}}= \dd_{h^{-1}}\ttt x= d(x)\dd_1
$$
\end{enumerate}
\end{lemma}
\begin{proof}
  Part (1) is immediate from the definition and
  $(\chi_i,\chi_j)=\de_{ij}$ (see \leref{l:bilinform3}); part (2) is
  trivial. Part (3) is proved in exactly the same way as in the proof
  of \leref{l:dim}.
  
  Parts (4), (5) are obtained from (2), (3) respectively by writing
  each side as a sum of homogeneus components. 
\end{proof}

\begin{lemma}\label{l:dim7}
  If $V$ is a simple module, and $\ph\colon V\to {}^gV, \psi\colon
  {}^gV\to V$,  then 
\begin{equation}
([\ph],[\psi])=\frac{1}{\dim V}\tr (\ph\psi).
\end{equation}
\end{lemma}
\begin{proof}
Consider $\ph\ttt \psi^*\colon V\ttt V^*\to {}^gV\ttt {}^gV^*$. Since
$V$ is simple, the multiplicity of $\1$ in $V\ttt V^*$ is one. Let
$i_V\colon \1\to V\ttt V^*$ be the canonical embedding; then  by
definition, the product $(x,y)$ is defined by 
$$
(\ph\ttt \psi^*)i_V=([\ph],[\psi])i_{{}^gV}.
$$
Pairing both sides with the canonical evaluation map $e_{{}^gV}\colon
{}^gV\ttt {}^gV^*\to\1$, we get
$$
d(V)([\ph],[\psi])=e_{{}^gV}(\ph\ttt \psi^*)i_V.
$$
It is immediate form the definition of the adjoint morphism that the
right-hand side of this identity is equal to $\tr(\ph\psi)$. 
\end{proof}
As before, this may be false if $V$ is not simple. 

\begin{lemma}\label{l:dim8}
For any $x,y\in \Ve$, one has
$$
(x,y)=d(x*y)=(\dd, x*y).
$$
\end{lemma}
\begin{proof}
Let $x\in \Ve_{g_1,h_1},
y\in \Ve_{g_2, h_2}$. If $g_1\ne g_2^{-1}$ or $h_1\ne h_2$, then both
sides are zero. Thus, it suffices to consider the case when $x=[\ph],
y=[\psi]$ for some $\ph\colon V_i\to {}^gV_i, \psi\colon {}^gV_j\to
{}V_j, V_i, V_j\in \CC_h$. In this case the result follows
from \leref{l:dim7}. 
\end{proof}

\subsection*{Extended Verlinde algebra and orbifold category}
Let us describe relation between the extended Verlinde
algebra $\Ve(\CC)$ and the Verlinde algebra $\V(\CG)$. Define the
map $F\colon \V(\CG)\to \Ve(\CC)$ by 
\begin{equation}\label{e:F}
F(X, \{\ph_g\})=\sum_g [\ph_g]. 
\end{equation}
Note that $[\ph_g]\in \Ve_{g^{-1}, *}$.

It is also useful to write  map $F$ in terms of the morphisms
$\ph_{g,h}\colon {}^hX\to {}^gX$ defined in \reref{r:ph}. 

\begin{lemma}\label{l:F3}
$$
F(X)=\frac{1}{|G|} \sum_{g,h}[\ph_{g,h}].
$$
\end{lemma}
\begin{proof}
It suffices to prove that for any $g\in G$, 
$$
[\ph_{a,b}]=[\ph_{ag,bg}]
$$
which in turn is equivalent to $[\ph_{a,b}]=[\ph_{1,ba^{-1}}]$. To
prove this, consider the following diagram: 
$$
\xymatrix{
{}^{ba^{-1}}X
     \ar[r]^{\ph_{1,ba^{-1}}}
     \ar[d]_{\ph_{b,ba^{-1}}}
 & X \ar[d]^{\ph_{a,1}}\\
{}^{b}X  \ar[r]^{\ph_{b,a}}
 & {}^aX
}
$$
By \eqref{e:ph3}, this diagram is commutative, and
$\ph_{b,ba^{-1}}=\Rho_b(\ph_{1,a^{-1}})$; thus, by \eqref{e:verrel22},
$[\ph_{a,b}=[\ph_{1,ba^{-1}}]$. 
\end{proof}

\begin{corollary}\label{c:F4}
$$
F(X)=\sum[\ph_g^{-1}]
$$
\end{corollary}

Also, define the map $G\colon \Ve(\CC)\to \V(\CG)$ as
follows: for $\psi\colon V\to {}^gV$, we let 
\begin{equation}\label{e:G}
\begin{aligned}
&G[\psi]=[f]\\
&f\colon \Ind V\to \Ind V\\
&f=\bigoplus_h R_h(\psi)\colon {}^hV\to  {}^{hg}V
\end{aligned}
\end{equation}

These definitions extend the maps $\V(\CG)\to \V(\CC)$, $\V(\CC)\to
\V(\CG)$ given by the restriction and induction functors defined by
\eqref{e:functors}.

\begin{example}\label{x:Gchig}
  By \exref{e:ind1}, $G[\1]=[\Ind\1]=[\F(G)\ttt \1]$. More generally,
  for the canonical morphism $\chi_0^g\colon\1\to {}^g\1$ defined by
  \eqref{e:chig}, one has $G\chi_0^g=\pi_g$, where $\pi_g\colon [\F(G)\ttt
  \1]\to [\F(G)\ttt \1]$ is the right regular action of $g$ on
  $\F(G)$. Decomposing $\F(G)=\bigoplus_\la \rho_\la\ttt \rho_\la^*$,
  where $\la$ runs over the set $\widehat G$ of isomorphism classes of
  irreducible $G$-modules, we see that
$$
G(\chi_0^g)=\bigoplus_\la \tr_{\rho_\la^*}(g)\cdot [\rho_\la\ttt \1].
$$

Since $\sum_g\pi_g$ acts by multiplication by $|G|$ on $\rho_0\ttt
\rho_0\subset \F(G)$ and by zero on $\rho_\la\ttt\rho_\la^*, \la\ne
0$, we get 
$$
G(\tilde\1)=G(\sum \chi_0^g)=|G|\cdot[\1_{\CG}]
$$
\end{example}

Relation between maps $F$ and $G$ is given by the following theorem,
parallel to \cite[Theorem 1.6]{KO}. 
\begin{theorem}\label{t:verlinde6}\indent\par 
\begin{enumerate}
\item 
$F(x\ttt y)=F(x)\ttt F(y)$, $G(F(x)\ttt y)=x\ttt G(y)$.
\item 
$GF(x)=|G| x$, $FG(x)=\sum_h R_h(x)$.
\item 
$F(x^*)=(F(x))^*, G(x^*)=(G(x))^*$, $[G(x)]_0=[x]_0$
\item 
$(Fx,y)=(x,Gy)$
\end{enumerate}
\end{theorem}
\begin{proof}
\begin{enumerate}
\item For $F$, immediate from the definition; for $G$, follows by a
  simple explicit computation. 
\item $FG(x)=\sum_h R_h(x)$ is immediate from the definitions;
  $GF(x)=|G| x$ follows from $GF(x)=G(F(x)\ttt \tilde\1)=x\ttt
  G(\tilde \1)$ (which is a special case of part (1)) and results of
  \exref{x:Gchig}.

\item For $*$, obvious; for $[\ ]_0$, note that if $\psi\colon V_i\to
  {}^gV_i$, $i\ne 0$, then the multiplicity $[\Res V_i:\1_{\CG}]=0$, so
  $[G\psi]_0=0=[\psi]_0$. For $\psi=\chi_0^g$, it follows from
  \exref{x:Gchig} that $[G\chi_0^g]_0=1=[\chi_0^g]$. 

\item Using results of previous parts, 
\begin{align*}
(Fx,y)&=[Fx\ttt y^*]_0=[G(Fx\ttt y^*)]_0=[x\ttt Gy^*]_0\\
&=(x,Gy)
\end{align*}
\end{enumerate}
\end{proof}

\begin{corollary}\label{c:verlinde7}
The map $F$ is an isomorphism 
$$
\V(\CG)\isoto (\Ve(\CC))^G
$$
\end{corollary}
\begin{proof}
  First, note that it is immediate from \eqref{e:ph3}and \leref{l:F3}
  that $\im F\subset (\Ve(\CC))^G$. The fact that it is an isomorphism
  is immediate from part (2) of \thref{t:verlinde6}: the inverse map
  is given by $\frac{1}{|G|}G$.
\end{proof}

\section{$S$-matrix}\label{s:smatrix}

Similar to the situation in $G=\{1\}$ case (see, e.g.,
\cite[Chapter 3]{BK}), in this section we introduce linear operators
$\tilde s, \tilde t\colon \Ve\to \Ve$. Later we will show that if
$\tilde s$ is non-degenerate, then after a simple renormalization these
operators define an action of the modular group $\SLZ$ on $\Ve$. 

In this section, we assume that $\CC$ is a $G$-equivariant fusion
category, and that the set $I$ of isomorphism classes of simple
objects in $\CC$ is finite. We denote by $\Ve$ the extended Verlinde
algebra of $\CC$ as defined in \seref{s:verlinde2}. 

Define
\begin{equation}\label{e:tmatrix}
\begin{aligned}
\tilde t\colon \Ve_{g,h}&\to \Ve_{gh,h}\\
          [\ph]&\mapsto[\th\ph]=[\ph\th]
\end{aligned} 
\end{equation}
where $\th$ is the universal twist. Equality $[\th\ph]=[\ph\th]$
follows from \eqref{e:verrel22} and $R_g\th=\th$ (\leref{l:twist}): if
$\ph\colon V\to {}^gV, V\in \CC_h$, then
\begin{equation}\label{e:thph=phth}
[\th_{{}^gV}\ph]=[V\xxto{\ph}{}^gV\xxto{\th}{}^{hg}V]
=[{}^{h^{-1}}V\xxto{\Rho_{h^{-1}}\th}V\xxto{\ph}{}^gV]
=[\ph\th_{{}^{h^{-1}}V}].
\end{equation}

\begin{figure}[th]
$$
(\tilde s[\ph])_i=d_i\quad
\raisebox{-0.5\height}{\begin{overpic}
{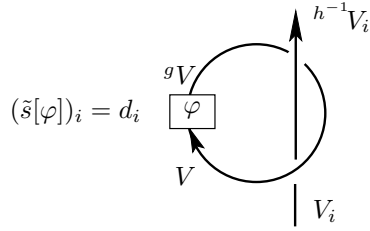}
\put(7,52){$\ph$}
\put(3,20){$V$}
\put(-2,65){${}^gV$}
\put(65,5){$V_i$}
\put(65,90){${}^{h^{-1}}\!V_i$}
\end{overpic}}
$$
\caption{$S$-matrix}\label{f:smatrix}
\end{figure}

We also define 
\begin{equation}\label{e:smatrix}
\begin{aligned}
\tilde s\colon\Ve_{g,h}&\to \Ve_{h^{-1}, g}\\
\tilde s[\ph]&=\sum_{i\in I_g} (\tilde s[\ph])_i
\end{aligned}
\end{equation}
where, for $\ph\colon V\to {}^g V, V\in \CC_h$, we define  $(\tilde
s[\ph])_i\colon V_i\to {}^{h^{-1}}V_i, V_i\in \CC_g$ by \firef{f:smatrix}.

Equivalently, the map $\tilde s$ can be described by the following
lemma. 
\begin{lemma}\label{l:sxy}
Let $\ph\colon V\to {}^gV, V\in \CC_h$ and $\psi\colon W\to {}^hW,
W\in \CC_g$, so $[\ph]\in\Ve_{g,h}, [\psi]\in\Ve_{h,g}$. Then  
$$
(\tilde s\ph,\psi)=
\raisebox{-0.5\height}{\begin{overpic}
{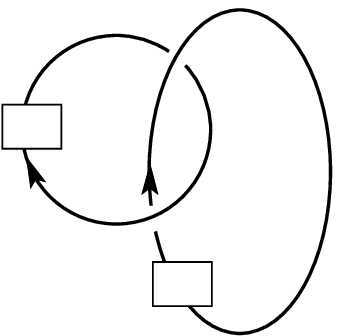}
\put(5,60){$\ph$}
\put(52,13){$\psi$}
\put(10,28){$V$}
\put(5,85){${}^gV$}
\put(85,40){$W$}
\put(53,25){${}^hW$}
\end{overpic}}
$$
\end{lemma}
\begin{proof}
Suffices to prove in the case when $W$ is simple, in which case it
follows from \leref{l:bilinform3}(4).
\end{proof}

\begin{lemma}\label{l:smatrix1}
\begin{enumerate}
\item The restriction of $\tilde s,\tilde t$ to $\Ve_{1,1}=\V(\CC_1)$
  coincides with the $\tilde s,\tilde t$ matrices defined in \cite{BK}.
\item $s,t$ are symmetric: $(\tilde s\ph,\psi)=(\ph, \tilde s\psi)$,
  $(\tilde t\ph, \psi)=(\ph,\tilde t\psi)$.
\item $\tilde s \1=\dd_1$, $\tilde s\tilde \1=\dd$. 
\item $\tilde s(x\ttt y)=\tilde s(y)*\tilde s(x)$

\end{enumerate}
\end{lemma}
\begin{proof}
\begin{enumerate}
\item  Immediate from the definition.

\item 
Symmetry of $\tilde s$ is immediate from
\leref{l:sxy}, after some simple manipulation of figures using the
Reidemeister moves of \seref{s:tangles}. To prove symmetry of $\tilde
t$, notice that by \leref{l:dim7}, one has 
$$
(\tilde tx,y)=d(\Theta*x*y)
$$
where  
$$
\Theta=\sum d_i [\theta_{V_i}]
$$
Thus, it suffices to prove that $\Theta$ is central with respect to
$*$, which is immediate from \eqref{e:thph=phth}. 

\item Immediate from the definition.
\item Let $\ph\colon V\to {}^gV, \psi\colon W\to {}^gW$. Using the
  graphical calculus of \seref{s:tangles}, we rewrite $\tilde s([\ph])*\tilde
  s([\psi])$ as follows:
\begin{align*}
&\tilde s([\ph])*\tilde s([\psi])=\\
&=\sum_i d_i \ 
\raisebox{-0.5\height}{\begin{overpic}
{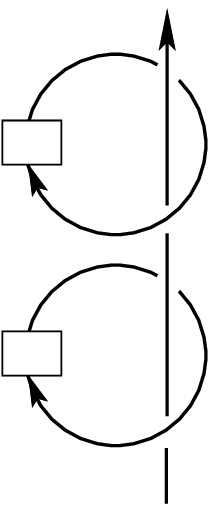}
\put(3,72){$\ph$}
\put(3,29){$\psi$}
\put(15,60){$V$}
\put(5,88){${}^gV$}
\put(15,17){$W$}
\put(3,46){${}^gW$}
\put(37,5){$V_i$}
\end{overpic}}
=\sum_i d_i \ 
\raisebox{-0.5\height}{\begin{overpic}
{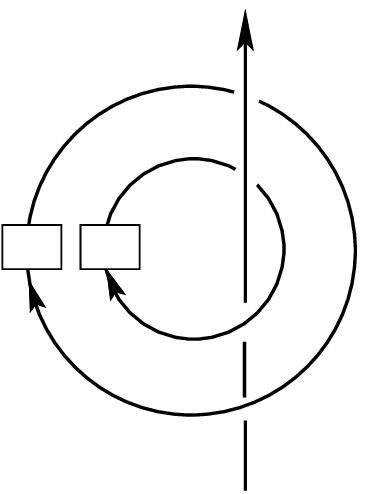}
\put(20,49){$\ph$}
\put(3,49){$\psi$}
\put(22,31){$V$}
\put(12,17){$W$}
\put(53,5){$V_i$}
\end{overpic}}
=\sum_i d_i\  
\raisebox{-0.5\height}{\begin{overpic}
{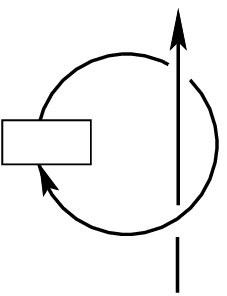}
\put(2,49){$\psi\!\ttt\!\ph$}
\put(12,10){$W\!\ttt\! V$}
\put(62,5){$V_i$}
\end{overpic}}\\
=&\tilde s([\psi]\ttt[\ph]).
\end{align*}

\end{enumerate}

\end{proof}

\begin{theorem}\label{t:smatrix2}
The operators $\tilde s, \tilde t$ satisfy the
  following relations 
\begin{equation}\label{e:modular1}
\begin{aligned}
&(\tilde s\tilde t)^3=p^+ \tilde s^2\\
&(\tilde s\tilde t^{-1})^3=p^- \tilde s^2 c\\
&c\tilde t=\tilde t c,\qquad c\tilde s=\tilde s c 
\end{aligned}
\end{equation}
where $c\colon \Ve\to \Ve$ is defined by $c[\ph]=[\ph^*]$ and 
\begin{equation}\label{e:p}
p^\pm=\sum_{i\in I_1} \th_i^{\pm 1}d_i^2.
\end{equation}
Note that this sum is over $I_1$, i.e. simple objects in $\CC_1$ only.
\end{theorem}

To porve this, theorem, we need first to prove several preparatory
results. 

\begin{proposition}\label{p:smatrix3}
For any $h\in G, i\in I_h$, 
$$
\raisebox{-0.5\height}{\begin{overpic}
{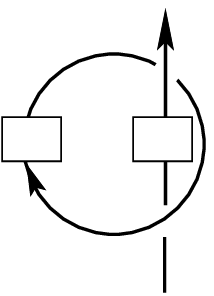}
\put(2,49){$\th^{\pm 1}$}
\put(48,49){$\th^{\pm 1}$}
\put(3,25){$\dd_h$}
\put(-15,74){${}^{h^{\pm 1}}\!\dd_h$}
\put(60,5){$V_i$}
\end{overpic}}\quad
=p^{\pm}\id_{V_i}
$$
where $\dd_h$ is as in \leref{l:dim6}. 
\end{proposition}
\begin{proof}
We prove the identity for $\th$; for $\th^{-1}$, the proof is similar. 

Since any morphism $V_i\to V_i$ is a multiple of identity, it suffices
to compute the trace of the left-hand side. Replacing $\dd_h$ by
$\dd_h^*=\dd_{h^{-1}}$ (see \leref{l:dim6}) and reversing the direction of the
corresponding strand, we see that the trace of left-hand side is given
by
$$
\raisebox{-0.5\height}{\begin{overpic}
{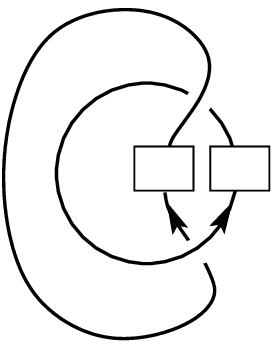}
\put(48,47){$\th$}
\put(68,47){$\th$}
\put(35,15){$\dd_{h^{-1}}$}
\put(65,5){$V_i$}
\end{overpic}}
$$
Using the identities $\th_{V\ttt W}=\th\ttt \th R^2$ (see
\leref{l:twist}) and $\dd_{h^{-1}}\ttt V_i=d_i \dd_1$ (see
\leref{l:dim6}), we can rewrite this as follows
$$
\raisebox{-0.5\height}{\begin{overpic}
{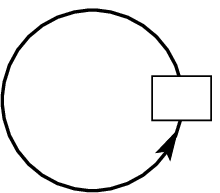}
\put(80,37){$\th$}
\put(20,-10){$\dd_{h^{-1}}\!\ttt V_i$}
\end{overpic}}\quad
= d_i\ 
\raisebox{-0.5\height}{\begin{overpic}
{smatrix4.eps}
\put(80,37){$\th$}
\put(80,10){$d_1$}
\end{overpic}}\quad
 =d_i\, p^+.
$$
\end{proof}

\begin{corollary}\label{c:smatrix4}
For any $V\in \CC_h$, one has 
$$
\raisebox{-0.5\height}{\begin{overpic}
{smatrix.eps}
\put(2,49){$\th^{\pm 1}$}
\put(3,20){$\dd_h$}
\put(65,5){$V$}
\put(65,90){${}^{h^{\mp 1}}\!V$}
\end{overpic}}\qquad
=p^\pm \th^{\mp 1}_V.
$$
\end{corollary}
\begin{proof}
Since both sides are functorial morphisms $V\to {}^{h^{\mp 1}}\!V$, it
suffices to prove this for $V$ being a simple module, in which case
it is immediate from \prref{p:smatrix3}.
\end{proof}

\begin{corollary}\label{c:smatrix5}
For any simple modules $V_i\in \CC_{h_1},  V_k\in\CC_{h_2}, h_1h_2=h$,
one has 
$$
\raisebox{-0.5\height}{\begin{overpic}
{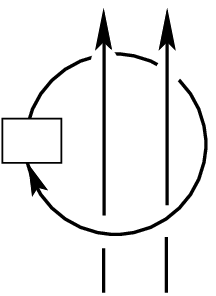}
\put(6,49){$\th$}
\put(3,20){$\dd_h$}
\put(38,5){$V_i$}
\put(58,5){$V_k$}
\end{overpic}}\qquad
=\quad
\raisebox{-0.5\height}{\begin{overpic}
{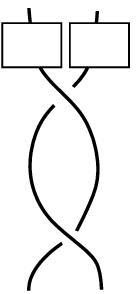}
\put(2,81){$\th^{-1}$}
\put(26,81){$\th^{-1}$}
\put(-2,5){$V_i$}
\put(35,5){$V_k$}
\end{overpic}}
$$
\end{corollary} 
Now we are ready to prove \thref{t:smatrix2}.

\begin{proof}[Proof of \thref{t:smatrix2}]
Let us prove the first identity, $(\tilde s\tilde t)^3=p^+ \tilde
s^2$. We rewrite it in the form 
$$
\tilde s \tilde t\tilde s= p^+ \tilde t^{-1}\tilde s\tilde t^{-1}.
$$

By definition, for $\ph \colon V_i\to {}^g V_i, V_i\in \CC_h$, we have
\begin{equation}\label{e:smatrix6}
\tilde s \tilde t\tilde s[\ph]=
\sum_{k\in I_{h^{-1}g}} d_k \sum_{j\in I_g} d_j\quad
\raisebox{-0.5\height}{\begin{overpic}
{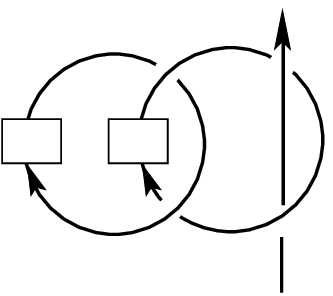}
\put(0,25){$V_i$}
\put(-4,62){${}^gV_i$}
\put(65,10){$V_j$}
\put(50,80){${}^{h^{-1}g}V_j$}
\put(87,10){$V_k$}
\put(90,85){${}^{g^{-1}}\!V_k$}
\put(5,45){$\ph$}
\put(40,44){$\th$}
\end{overpic}}
\end{equation}

Using Reidemeister moves for framed graphs, we rewrite
(\ref{e:smatrix6}) in the form 
 \begin{equation}\label{e:smatrix7}
\tilde s \tilde t\tilde s[\ph]=
\sum_{k\in I_{h^{-1}g}} d_k \sum_{j\in I_g} d_j\quad
\raisebox{-0.5\height}{\begin{overpic}
{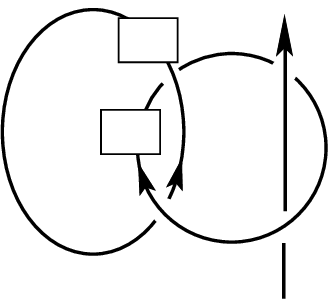}
\put(60,50){$V_i$}
\put(-4,80){${}^gV_i$}
\put(65,8){$V_j$}
\put(87,10){$V_k$}
\put(90,85){${}^{g^{-1}}V_k$}
\put(41,78){$\ph$}
\put(38,48){$\th$}
\end{overpic}}
\end{equation}

Using \coref{c:smatrix5}, this can be rewritten as 

 \begin{equation}\label{e:smatrix8}
\tilde s \tilde t\tilde s[\ph]=p^+
\sum_{k\in I_{h^{-1}g}} d_k\quad 
\raisebox{-0.5\height}{\begin{overpic}
{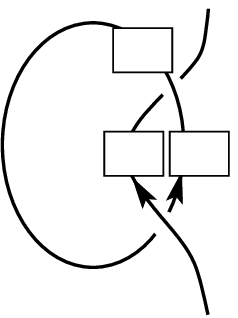}
\put(60,62){$V_i$}
\put(-5,20){${}^gV_i$}
\put(68,10){$V_k$}
\put(43,83){$\ph$}
\put(34,48){$\th^{-1}$}
\put(55,48){$\th^{-1}$}
\end{overpic}}
\quad=p^+\quad 
\raisebox{-0.5\height}{\begin{overpic}
{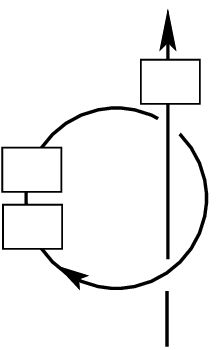}
\put(5,18){$V_i$}
\put(-10,66){${}^{h^{-1}g}V_i$}
\put(50,10){$V_k$}
\put(52,90){${}^{g^{-1}}V_k$}
\put(5,34){$\ph$}
\put(2,48){$\th^{-1}$}
\put(42,74){$\th^{-1}$}
\end{overpic}}
\quad=p^+ (\tilde t^{-1}\tilde s\tilde t^{-1})[\ph].
\end{equation}

This proves the first identity of \thref{t:smatrix2}. The second is
proved similarly, using \coref{c:smatrix4} for $\th^{-1}$. 

\end{proof}

Finally, the following theorem establishes a relation between the
$\tilde s$ operator for the $G$-equivariant category $\CC$ and the
orbifold category $\CG$. The following theorem is an analog of
\cite[Theorem 4.1]{KO}; however, the use of extended Verlinde algebra
allows us to simplify the statement of this theorem. 

\begin{theorem}\label{t:smatrix9}
Let $F\colon \V(\CG)\to\Ve (\CC), G\colon \Ve(\CC)\to \V(\CG)$
be defined by \eqref{e:F}, \eqref{e:G}. Then 
\begin{alignat*}{2}
F\tilde s&= |G|\cdot \tilde s F &\qquad F\tilde t&=\tilde t F\\
G\tilde s&= \frac{1}{|G|}\tilde s G &\qquad G\tilde t&=\tilde t G
\end{alignat*}
\end{theorem}
\begin{proof}
  
  We start by proving commutation relations of $F,G$ with $\tilde t$.
 To prove $F\tilde t=\tilde t F$, recall that by definition, for
 $x=[(X,\{\ph_g\})]\in \V(\CG)$, one has 
$$
F(\tilde t x)=\sum_g[\th^{\CG}\ph_g]
$$
where $\th^{\CG}$ is the twist in $\CG$. Using definition of
twist in $\CG$ (\thref{t:orbicategory2}), we rewrite it as
$$
F(\tilde t x)=
\sum_{g}[{}^gX\xxto{\ph_g}X\xxto{\th}{}^*X\xxto{\ph_*}X]
$$
where $*$ has the following meaning: if we write this sequence as a
direct sum of homogeneous sequences (i.e., taking place in $\CC_h,
h\in G$), then on $\CC_h$, $*=h$. 

Using functoriality of $\th$ and definition of $\CG$, we rewrite it
as follows: 
\begin{align*}
F(\tilde t x)=&
\sum_g    [{}^gX\xxto{\th}{}^{*g}X\xxto{\Rho_*\ph_g}{}^*X\xxto{\ph_*}X]\\
=&\sum_g    [{}^gX\xxto{\th}{}^{*g}X\xxto{\ph_{*g}}X]
=\sum_a    [{}^{*^{-1}a}X\xxto{\th}{}^{a}X\xxto{\ph_a}X]\\
=&\sum_a [\ph_a\th]
=\tilde t (F x).
\end{align*}

To prove identity $G\tilde t=\tilde t G$, recall that by definition, 
for $\psi\colon V\to {}^gV, V\in \CC_h$, we have 
$$
G(\tilde t [\psi])=G([\th\psi])=[f]
$$
where $f\colon \bigoplus_a {}^aV\to \bigoplus_a {}^aV$ is defined by 
\begin{align*}
&\bigoplus_a {}^aV\xxto{\oplus\Rho_a(\th\psi)}
  \bigoplus_a {}^{ahg}V\xxto{permutation}\bigoplus_a {}^aV\\
&\qquad=\bigoplus_a {}^aV\xxto{\oplus\Rho_a(\psi)}
  \bigoplus_a {}^{ag}V
\xxto{\th}
  \bigoplus_a {}^{ahg}V
\xxto{permutation}\bigoplus_a {}^aV
\end{align*}
On the other hand, $\tilde t G[\psi]$ is by definition given by the
class of the following morphism  
\begin{align*}
&\bigoplus_a {}^aV\xxto{\oplus\Rho_a(\psi)}
  \bigoplus_a {}^{ag}V\xxto{\th^{\CG}}\bigoplus_a {}^{ag}V
\xxto{permutation}\bigoplus_a {}^aV\\
&\qquad=\bigoplus_a {}^aV\xxto{\oplus\Rho_a(\psi)}
  \bigoplus_a {}^{ag}V
\xxto{\th}
  \bigoplus_a {}^{ahg}V
\xxto{permutation}\bigoplus_a {}^aV
\end{align*}
which proves $G\tilde t[\psi]=\tilde t G[\psi]$.

To prove
identities involving $\tilde s$, it suffices to prove that for any
$x\in \V(\CG), y\in \Ve(\CC)$ one has
\begin{equation}\label{e:smatrix9}
(\tilde s Fx,y)_\CC=\frac{1}{|G|}(\tilde s x, Gy)_{\CG}
\end{equation}

Indeed, since $F$ and $G$ are adjoint (see \thref{t:verlinde6}),
\eqref{e:smatrix9} implies $(\tilde s Fx,y)=\frac{1}{|G|}(F\tilde
sx,y)$, and thus, since the form is non-degenerate, 
$\tilde s Fx=\frac{1}{|G|}F\tilde s x$. Similarly, using  by symmetry of
$\tilde s$ (see \leref{l:smatrix1}) and adjointness of $F,G$, we see that
left-hand side of    \eqref{e:smatrix9} is equal to $(Fx, \tilde
sy)=(x, G \tilde s y)$, and the right-hand side is equal to $
\frac{1}{|G|}(x, \tilde s Gy)$, so \eqref{e:smatrix9} implies 
$G \tilde s y=\frac{1}{|G|}\tilde s Gy$.

So let us prove \eqref{e:smatrix9}. Without loss of generality, we can
assume that $x=[X,\{\ph_g\}], (X,\{\ph_g\})\in \CG$, and $y=[\psi]\in
\Ve_{h,g}$, where $\psi\colon W\to {}^hW, W\in \CC_g$, $hg=gh$.  In
this case, using \leref{l:sxy} and formula for  $F$ given in 
\coref{c:F4}, we see that the left-hand side of \eqref{e:smatrix9} is
given by
\begin{equation}\label{e:smatrix10}
(\tilde s Fx,y)=
\raisebox{-0.5\height}{\begin{overpic}
{smatrix10.eps}
\put(1,61){$\ph_g^{-1}$}
\put(52,13){$\psi$}
\put(10,28){$X_h$}
\put(1,88){${}^gX_h$}
\put(85,40){$W$}
\put(53,25){${}^hW$}
\end{overpic}}
\end{equation}
where $X_h$ is the component of $X$ in $\CC_h$ (all other components
give zero contribution). 

Now let us compute the right-hand side. Recall that by definition,
$G[\psi]=[f]$, where $f\colon \Ind W\to \Ind W$ is defined by
\eqref{e:G}. Using an analog of \leref{l:sxy} for usual fusion
categories (which can be obtained from \leref{l:sxy} by letting
$G=\{1\}$), we see that
$$
(\tilde s x, Gy)= 
\raisebox{-0.5\height}{\begin{overpic}
{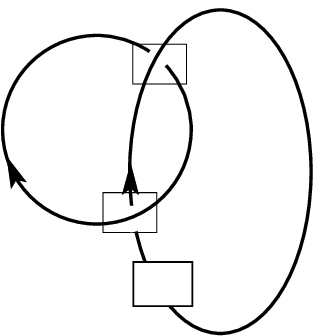}
\put(47,13){$f$}
\put(10,28){$X$}
\put(100,40){$\Ind W$}
\end{overpic}}
$$
where the boxed crossings are commutativity isomorphisms in $\CG$. 

Writing $X=\bigoplus X_a, X_a\in \CC_a$ and $\Ind W=\bigoplus {}^bW$
and using definition of the commutativity isomorphism in $\CG$, we
can rewrite it as follows: 
\medskip

$$
(\tilde s x, Gy)= \sum_{a,b: a^{-1}bh=b}
\raisebox{-0.5\height}{\begin{overpic}
{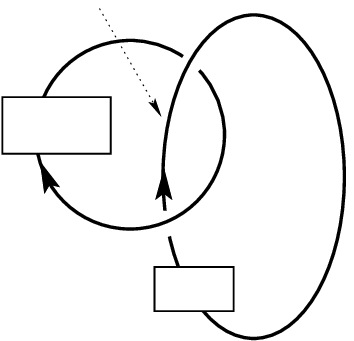}
\put(47,13){$R_b\psi$}
\put(10,28){$X_a$}
\put(2,60){$\ph^{-1}_{bgb^{-1}}$}
\put(100,40){${}^bW$}
\put(17,100){${}^{a^{-1}bh}W$}
\put(53,25){${}^{bh}W$}
\end{overpic}}
$$
(as before, one easily sees that components with $a\ne bhb^{-1}$ give
zero contribution). Replacing in this formula  $a$ by $bhb^{-1}$, we
get 
\medskip
\begin{align*}
(\tilde s x, Gy)= &\sum_{b} 
\raisebox{-0.5\height}{\begin{overpic}
{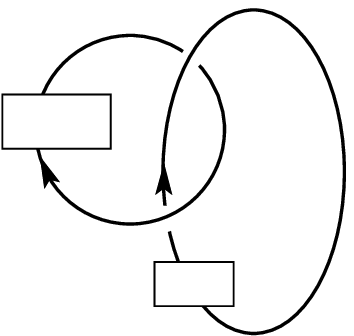} 
\put(47,13){$R_b\psi$} 
\put(10,28){$X_{bhb^{-1}}$} 
\put(2,60){$\ph^{-1}_{bgb^{-1}}$} 
\put(100,40){${}^bW$} 
\put(53,25){${}^{bh}W$} 
\end{overpic}}\qquad 
=
\sum_{b}R_b\Biggl(\quad
\raisebox{-0.5\height}{\begin{overpic}
{smatrix16.eps} 
\put(51,13){$\psi$} 
\put(10,28){$X_h$} 
\put(0,80){${}^gX_h$} 
\put(6,60){$\ph_{h}^{-1}$}
\put(100,40){$W$} 
\put(53,25){${}^{h}W$}
\end{overpic}}
\quad\Biggr)\\ 
&=|G|\quad
\raisebox{-0.5\height}{\begin{overpic}
{smatrix16.eps} 
\put(51,13){$\psi$}
\put(10,28){$X_h$}
\put(0,80){${}^gX_h$} 
\put(6,60){$\ph_{h}^{-1}$}
\put(100,40){$W$}
\put(53,25){${}^{h}W$}
\end{overpic}}
\end{align*}
Comparing this with \eqref{e:smatrix10}, we get \eqref{e:smatrix9}.
\end{proof}

\section{Modular equivariant categories}\label{s:modular} 

Throughout this section, $\CC$ is a $G$-equivariant fusion category
with finitely many isomorphism classes of simple objects, and
$\Ve=\Ve(\CC)$ is the extended Verlinde algebra as defined in
\seref{s:verlinde2}.

\begin{definition}\label{d:modular}
  A $G$-equivariant fusion category with finitely many isomorphism
  classes of simple objects is called {\em modular} if the
  operator $\tilde s\colon \Ve\to\Ve$, defined by \eqref{e:smatrix},
  is invertible.
\end{definition}
This definition generalizes the well-known definition of a modular
tensor category

\begin{remark}
  This definition is different from the one given in \cite{T2}.
  Namely, the definition of \cite{T2} only requires that the
  subcategory $\CC_1$ be modular. It is easy to see (see \thref{t:modular3}
  below and discussion following it) that modularity of $\CC$ implies
  modularity of $\CC_1$ but converse is not true. Thus, our definition
  is stronger than that of \cite{T2}.
\end{remark}

\begin{theorem}\label{t:modular2}
Let $\CC$ be a $G$-equivariant modular category. Then:
\begin{enumerate}
\item The numbers $p^\pm$ defined by \eqref{e:p} are non-zero.
\item Let 
\begin{equation}
s=D^{-1}\tilde s,\qquad t=\zeta^{-1}\tilde t
\end{equation}
where $\tilde s,\tilde t\colon \Ve\to \Ve$ are defined
in \seref{s:smatrix} and  
\begin{equation}\label{e:Dzeta}
D=\sqrt{p^+p^-},\qquad \zeta=(p^+p^-)^{1/6}.
\end{equation}

Then so defined $s,t$ satisfy the relations of $\SLZ$: 
\begin{equation}\label{e:slz}
(st)^3=s^2,\quad s^2=c,\quad ct=tc,\quad c^2=1
\end{equation}
where $c\colon \Ve\to\Ve$ is as in \thref{t:smatrix2}.
\end{enumerate}
\end{theorem}
\begin{proof}
Rewriting equalities of  \thref{t:smatrix2} in the form
\begin{align*}
&\tilde s\tilde t \tilde s =p^+\tilde t^{-1}\tilde s \tilde t^{-1}\\
&\tilde s\tilde t^{-1} \tilde s =p^- c \tilde t\tilde s \tilde t
=p^- c \tilde t \cdot \tilde s \tilde t \tilde s \cdot \tilde s^{-1}
\end{align*}
and substituting the first equality into the second one, we get
$\tilde s^2=p^+p^- c$. After this, the results immediately follow from
\thref{t:smatrix2}.
\end{proof}

Note that the numbers $\zeta, D$ are the same as for the modular
category $\CC_1$; in other words, the central charge and rank of a
$G$-equivariant modular category are the same as for its neutral
part.

As one might expect, modularity of $\CC$ is closely related with
modularity of the orbifold category $\CG$ and with the modularity of
the untwisted sector $\CC_1$.

\begin{theorem}\label{t:modular3}
Let $\CC$ be a $G$-equivariant modular category. Then both $\CG$
and $\CC_1$ are modular categories, and the numbers $D,\zeta$ defined
by \eqref{e:Dzeta} are related by 
\begin{align*}
&\zeta(\CC_1)=\zeta(\CC)=\zeta(\CG)\\
&D(\CC_1)=D(\CC)=\frac{D(\CG)}{|G|}
\end{align*}
\end{theorem}

\begin{proof}
Modularity of $\CC_1$ is immediate, as the $\tilde s$-matrix for
$\CC_1$ is just the restriction of $\tilde s$-operator for $\CC$ to
$\V(\CC_1)=\Ve_{1,1}$. 

To prove modularity of $\CG$, note that by \coref{c:verlinde7}
and \thref{t:smatrix9},  we have an embedding
$F\colon \V(\CG)\injto \Ve(\CC)$ and $F\tilde s=|G|\cdot \tilde s
F$. Thus, if $\tilde s x=0$ for some $x\in \V(\CG)$, then $Fx\in
\Ker \tilde s$ in $\Ve(\CC)$, which contradicts modularity of $\CC$. 

To prove the relation between numbers $D,\zeta$ for $\CC, \CG$,
note that applying $F$ to both sides of the identity $(\tilde s\tilde
t)^3=p^+\tilde s^2$ for $\CG$ and using \thref{t:smatrix9}, we get
$p^+(\CC)=|G|\cdot  p^+(\CG)$; in a similar way we get
$p^-(\CC)=|G|\cdot  p^-(\CG)$. 
\end{proof}    
A natural question is whether converse is also true, i.e. whether
modularity of $\CC_1$ or $\CG$ imply modularity of $\CC$. One easiy
sees that modularity of $\CC_1$ does not imply modularity of $\CC$:
for example, if one takes any (usual) modular category $\CC$ and
considers it as a $G$-equivariant fusion category with trivial
$G$-grading and trivial action of $G$, then $\CC_1$ is modular, but
$\CC$ is not a $G$-equivariant modular category (which can be easily
seen from \leref{l:modular4} below).

On the other hand, modularity of $\CG$ does imply modularity of
$\CC$.

\begin{theorem}\label{t:modular4}
A $G$-equivariant fusion category $\CC$ is modular iff the orbifold
category $\CG$ is modular. 
\end{theorem}
\begin{proof}
  By \thref{t:modular3}, if $\CC$ is modular then so is $\CG$. Thus,
  we only need to prove that if $\CG$ is modular, then $\CC$ is
  modular.

Assume that $\CG$ is modular; then  $\tilde s^2_{\CG}=D^2 c$ for some
non-zero $D=D(\CG)$. Applying to this $F$, we get 
$$
F\tilde s^2_{\CG}=D^2cF=\frac{1}{|G|^2}\tilde s^2_{\CC}F.
$$
Applying this to $\1\in \V(\CG)$ and using $F\1=\tilde \1$, we get 
\begin{equation}\label{e:s^2}
\tilde s^2_{\CC}\tilde\1=\mu\tilde \1,\qquad \mu=|G|^2D^2\ne 0
\end{equation} 
(compare with \cite[Lemma 4.6]{KO}). 

Let us now show that for any $x\in \Ve(\CC), \tilde s x\ne 0$. Indeed,
let us compute $(\tilde s^2\tilde\1, x\ttt y^*)$. On one hand, by
\eqref{e:s^2} and definition of the bilinear form, we get
 $$
(\tilde s^2\tilde\1, x\ttt y^*)=\mu (\tilde \1, x\ttt y^*)=\mu(x,y).
$$
On the other hand, using \leref{l:smatrix1}, we get
$$
(\tilde s^2\tilde\1, x\ttt y^*)
=(\tilde s\tilde\1, \tilde s(x\ttt y^*))
=(\tilde s\tilde\1, \tilde s(y^*)*\tilde s(x))
$$

Thus, if $\tilde s(x)=0$, then $(\tilde s^2\tilde\1, x\ttt
y^*)=\mu(x,y)=0$ for all $y$, which contradicts non-degeneracy of the
form $(x,y)$ (\leref{l:bilinform3}). 
\end{proof}
\begin{remark}
Combining \thref{t:modular3}, \thref{t:modular4}, we get that if
$\CG$ modular, then $\CC_1$ is also modular, which is exactly a
statement of \cite[Theorem 4.5]{KO} in our situation. In fact, the
proof given above is parallel to the proof in \cite{KO}; we used the
bilinear form to simplify some of the arguments in \cite{KO}. 
\end{remark}

As in the usual case, modularity implies a number of remarkable
properties of the category. Here are the most immediate ones. 

\begin{lemma}\label{l:modular4}
Let $\CC$ be a  modular $G$-equivariant category. Then 
$$
|I_g|=|(I_1)^g|=\bigl|\{i\in I_1\st V_i\simeq
{}^gV_i\}\bigr|. 
$$
In particular, for every $g$, $|I_g|\ge 1$. 
\end{lemma}
\begin{proof}
  In a modular $G$-equivariant category, $s$-matrix gives an
  isomorphism $\Ve_{g,h}\isoto \Ve_{h^{-1},g}$. Thus, 
$$
|I_g|=\dim \Ve_{1,g}=\dim\Ve_{g,1}=\bigl|\{i\in I_1\st V_i\simeq
{}^gV_i\}\bigr|. 
$$
\end{proof}

\begin{lemma}\label{l:duality}
For any $x,y\in \Ve$, one has 
\begin{align*}
s(x\ttt y)&= D\, s(y)* s(x)\\
s(x*y)=\frac{1}{D}s(x)\ttt s(y)
\end{align*}
\end{lemma}
\begin{proof}
The first identity is immediate from \leref{l:smatrix1} and definition
of $s$. To prove the second identity, write $x=s(a), y=s(b)$ for some
$x,y$ (which is possible because $s$ is invertible); then, using the
first identity, we get  
\begin{align*}
s(x*y)&=s(s(a)*s(b))=\frac{1}{D}s^2(b\ttt a)=\frac{1}{D}s^2(a)\ttt
s^2(b)\\
&=\frac{1}{D}s(x)\ttt s(y)
\end{align*}
\end{proof}

\begin{corollary}\label{c:duality}
Let $\CC$ be a  modular $G$-equivariant category. Then the two
structures of associative algebra on $\Ve$, defined by $\ttt$ and
$*$, are isomorphic.
\end{corollary}

\begin{lemma}\label{l:holomorphic}
If $\CC_1=\Vec$, then $\CC$ is equivalent as a
monoidal  category to the category $\GVec^\om$ of 
twisted $G$-graded vector spaces for some $\om\in H^3(G, \Ctimes)$.  
\end{lemma}
\begin{proof}
  It follows from \leref{l:modular4} that $|I_g|=1$ for all $g$; in
  other words for every $g$, there is a unique simple object $X_g\in
  \CC_g$. Since the dual of a simple object is simple, this implies
  $X_g^*\simeq X_{g^{-1}}$ (non-canonically).

  Consider $X_g\ttt X_g^*$. This object lies in $\CC_1$, thus it must
  be a multiple of $\1$. On the other hand, for a simple object $X$,
  multiplicity of $\1$ in $X\ttt X^*$ is one. Thus, $X_g\ttt
  X_g^*=\1$. 
 
  Consider tensor product $X_g\ttt X_h$. Since $X_g\ttt X_h\in
  \CC_{gh}$, we have $X_g\ttt X_h \simeq N_{gh}X_{gh}$ for
  some multiplicities $N_{gh}\in \Z_+$. Tensoring both sides with $X_h^*$ and
  using $X_h\ttt X_h^*\simeq \1$, we get $X_g \simeq N_{gh}X_{gh}\ttt
  X_h^*$. Since $X_g$ is simple, this implies $N_{gh}=1$,
  i.e. $X_g\ttt X_h\simeq X_{gh}$. 
\end{proof}

This lemma, combined with \exref{x:twistedvecG}, gives a much
simpler proof of the main result of \cite{orbi1} (see \cite[Corollary
5.13]{orbi1}). The only difference is that here, we used the
assumption that $\CC$ is modular, whereas \cite{orbi1} used the
assumption that $\CG$ is modular.

\begin{lemma}\label{l:verlindef}
For $x\in \Ve$, let $L_x$ be the operator of left multiplication
by $x$: $L_x y=x\ttt y$. Then 
$$
sL_xs^{-1}=D_x
$$
where $D_x$ is the renormalized operator of $*$-multiplication by $s(x)$:
$$
D_x y = y*\frac{s(x)}{D} 
$$
\end{lemma}
\begin{proof}
This is equivalent to 
$$
s(x\ttt s^{-1}(y))=\frac{y*s(x)}{D}
$$
which immediately follows from \leref{l:duality}. 
\end{proof}

This lemma is an analog of the famous Verlinde formula for the
usual modular categories.  Note, however, that in $G$-equivariant
case, both $\ttt$ and $*$ are in general non-commutative, and operator
$D_x$ is not diagonal; thus, we can not say that ``$s$-matrix
diagonalizes fusion rules''. However, in some special cases (namely,
when $G$ is commutative and the cohomology classes $\al$ in
\thref{t:algebra3} are trivial), $*$ is commutative and therefore we
can use \leref{l:verlindef} to compute the fusion coefficients in
$\Ve$. This will be described in a forthcoming paper.



\end{document}